# HIGHER ORDER SEMIPARAMETRIC FREQUENTIST INFERENCE WITH THE PROFILE SAMPLER[1]

By Guang Cheng and Michael R. Kosorok

*Duke University and University of North Carolina at Chapel Hill*

We consider higher order frequentist inference for the parametric component of a semiparametric model based on sampling from the posterior profile distribution. The first order validity of this procedure established by Lee, Kosorok and Fine in [*J. American Statist. Assoc.* **100** (2005) 960–969] is extended to second-order validity in the setting where the infinite-dimensional nuisance parameter achieves the parametric rate. Specifically, we obtain higher order estimates of the maximum profile likelihood estimator and of the efficient Fisher information. Moreover, we prove that an exact frequentist confidence interval for the parametric component at level $\alpha$ can be estimated by the $\alpha$-level credible set from the profile sampler with an error of order $O_P(n^{-1})$. Simulation studies are used to assess second-order asymptotic validity of the profile sampler. As far as we are aware, these are the first higher order accuracy results for semiparametric frequentist inference.

**1. Introduction.** The focus of this paper is on higher order frequentist inference for the parametric component $\theta$ of a semiparametric model. In addition to the $d$-dimensional Euclidean parameter $\theta$, semiparametric models also have an infinite-dimensional parameter $\eta$, sometimes called the "nuisance" parameter. A classic example is the Cox proportional hazards model for right-censored survival data [4], where interest is focused on the log hazard ratios $\theta$ for the regression covariate vector $z$. The integrated baseline hazard function $\eta$ is the infinite-dimensional nuisance parameter. The involvement of an infinite-dimensional nuisance parameter in semiparametric models generally complicates maximum likelihood inference for $\theta$. In partic-

Received June 2007; revised June 2007.
[1]Supported in part by Grant CA075142.
*AMS 2000 subject classifications.* Primary 62G20, 62F25; secondary 62F15, 62F12.
*Key words and phrases.* Higher order frequentist inference, posterior distribution, Markov chain Monte Carlo, profile likelihood, Cox proportional hazards model, proportional odds model, case-control studies with a missing covariate.

ular, estimating the limiting variance of $\sqrt{n}(\hat{\theta}_n - \theta_0)$, where $\theta_0$ is the true value of $\theta$, usually requires estimating an infinite-dimensional operator.

The related studies concerning higher order frequentist inference in the parametric models under the Bayesian set-up focus on the choice of priors, such as objective priors [30]. However, it turns out that extending the objective prior approach to the semiparametric setting seems to require a higher-than-second order expansion of the profile likelihood and appears to be quite difficult. A similar hurdle appears to be required for extending the higher order bootstrap results for parametric models [8] to the semiparametric setting. Interestingly, general first-order bootstrap results for semiparametric M-estimators have only recently been developed (see [31] and [18]). Higher order extensions for any of these approaches would be very useful. However, in this paper, we will pursue an apparently simpler approach to obtaining higher order likelihood inference for semiparametric models based on the profile sampler proposed in [17].

The profile sampler provides a first-order correct approximation of the maximum likelihood estimator $\hat{\theta}_n$ and consistent estimation of the efficient Fisher information for $\theta$ based on sampling from the posterior of the profile likelihood, even when the nuisance parameter is not estimable at the $\sqrt{n}$ rate. The validity of the profile sampler relies on special properties of the profile likelihood in semiparametric models, some of which are extensively studied in [20, 21] and [22]. The profile likelihood for the parameter $\theta$ is $pl_n(\theta) = \sup_\eta \mathrm{lik}_n(\theta, \eta)$, where $\mathrm{lik}_n(\theta, \eta)$ is the full likelihood given $n$ observations. We also define $\hat{\eta}_\theta = \arg\max_\eta \mathrm{lik}_n(\theta, \eta)$. The maximum likelihood estimator for the full likelihood is thus $(\hat{\theta}_n, \hat{\eta}_n)$, where $\hat{\eta}_n = \hat{\eta}_{\hat{\theta}_n}$. Consideration of the profile likelihood in frequentist inference about $\theta$ can be traced back to the ordinary parametric model. An intuitive interpretation for the validity of the profile likelihood in semiparametric models is that it can be viewed as an estimator of the least favorable submodel for the estimation of $\theta$ [25]. The least favorable submodel, which will be briefly introduced in the next section, is the closest parametric model to the semiparametric model in the sense of information. In practice, the profile likelihood can often be easily computed using procedures such as the stationary point algorithm (as used in, e.g., [14]) or the iterative convex minorant algorithm introduced in [7] to find $\hat{\eta}_\theta$ if $\eta$ is a monotone function.

An advantage of the profile sampler is that a prior on the infinite-dimensional parameter is not required to obtain valid frequentist inference about $\theta$. Assigning a prior on $\eta$ can be quite challenging since for some models, there is no direct extension of the concept of a Lebesgue dominating measure for the infinite-dimensional parameter set involved [15]. The fully Bayesian approach can obviously be the basis for inference on $\theta$ alone via the marginal posterior. The first-order valid results in [26] indicate that the marginal semiparametric posterior is asymptotically normal and centered at the corresponding maximum likelihood estimator or posterior mean, with covariance matrix equal to the inverse of the efficient Fisher information. Unfortunately, this marginal approach does not circumvent the need to specify a prior on $\eta$, with all of the difficulties that entails.

The main contribution of this paper is the development of higher order frequentist inference for the parametric component of a semiparametric model through the profile sampler procedure proposed in [17] by assuming stronger assumptions on the semiparametric model and prior. We assume that the nuisance parameters of the semiparametric models studied in this paper have the parametric rate. This assumption permits the treatment of the likelihood as essentially parametric in certain aspects. This enables the second-order frequentist inference results for parametric models to be naturally extended to the semiparametric setting, although we note that considerable technical difficulties are present despite this simplification. To accomplish the above higher order inference, we require stricter—but still reasonable—regularity conditions than those imposed by [22] on the least favorable submodel. This is reviewed in Section 2. The initial technical step, presented in Section 3, is to establish higher order versions of expansions (5)–(6) in [22]. In Section 4, we find that the mean (median) value and the inverse of the variance of the MCMC chain from the profile sampler are actually higher order estimates of the maximum likelihood estimator and the efficient Fisher information, respectively. The main result of Section 4 is to prove that an exact frequentist confidence interval for $\theta_0$ can be estimated by the credible set from the profile sampler with an error of order only $O_P(n^{-1})$. In Section 5, we discuss three examples and some simulation results. Section 5 is followed by a discussion in Section 6 of future research interests. We postpone most of the technical details to the proofs given in Section 7.

As far as we are aware, these are the first higher order accuracy results for semiparametric frequentist inference. This is quite distinct from the concept of second-order efficiency in semiparametric models (see [9] and [5]) which we do not consider in this paper. The two important tools we use in this paper are separately empirical processes and sandwich techniques [22], with which we establish upper and lower bounds for the error in the profile log-likelihood expansion. For ease of exposition, we assume throughout the paper that $\theta \in \mathbb{R}^1$. However, the results can be readily extended to higher dimensions. The confidence "interval" and credible set for $d$-dimensional $\theta$ are a rectangle, a cuboid and a hypercuboid when $d = 2$, $d = 3$ and $d \geq 4$, respectively.

**2. Preliminaries.** We assume the data $X_1, \ldots, X_n$ are i.i.d. throughout the paper. The sample space $\mathcal{X}$ will depend on the semiparametric model which is defined by a density $\{p_{\theta,\eta}(x) : \theta \in \Theta, \eta \in \mathcal{H}\}$, where $\mathcal{H}$ is an arbitrary subset that will typically be infinite-dimensional. We first review the concept of a least favorable submodel and then present some notation and assumptions that will be used throughout the paper.

2.1. *The least favorable submodel.* The *score function* for $\theta$, $\dot{\ell}_{\theta,\eta}$ is defined as the partial derivative with respect to $\theta$ of the log-likelihood given fixed $\eta$. A score function for $\eta_0$ is of the form

$$\left.\frac{\partial}{\partial t}\right|_{t=0} \log p_{\theta_0,\eta_t}(x) \equiv A_{\theta_0,\eta_0} h(x),$$

where $h$ is a "direction" by which $\eta_t \in \mathcal{H}$ approaches $\eta_0$, running through some index set $H$. $A_{\theta,\eta}: H \mapsto L_2^0(P_{\theta,\eta})$ is the score operator for $\eta$. The *efficient score function* for $\theta$ is defined as $\tilde{\ell}_{\theta,\eta} = \dot{\ell}_{\theta,\eta} - \Pi_{\theta,\eta}\dot{\ell}_{\theta,\eta}$, where $\Pi_{\theta,\eta}\dot{\ell}_{\theta,\eta}$ minimizes the squared distance $P_{\theta,\eta}(\dot{\ell}_{\theta,\eta} - k)^2$ over all functions $k$ in the closed linear space of the score functions for $\eta$ (the "nuisance scores"). The variance of $\tilde{\ell}_{\theta,\eta}$, called the *efficient information matrix*, $\tilde{I}_{\theta,\eta}$, is the Cramér–Rao bound for estimating $\theta$ in the presence of the infinite-dimensional nuisance parameter $\eta$. We denote $\tilde{\ell}_{\theta_0,\eta_0}$ and $\tilde{I}_{\theta_0,\eta_0}$ by $\tilde{\ell}_0$ and $\tilde{I}_0$, respectively.

A submodel $t \mapsto p_{t,\eta_t}$ is defined to be *least favorable* at $(\theta, \eta)$ if $\tilde{\ell}_{\theta,\eta} = \partial/\partial t \log p_{t,\eta_t}$, given $t = \theta$. The "direction" along which $\eta_t$ approaches $\eta$ in the least favorable submodel is called the *least favorable direction*. Generally, the least favorable direction at $(\theta, \eta)$ in semiparametric models can be obtained by solving for $h_{\theta,\eta}$ in the equation $P(\dot{\ell}_{\theta,\eta} - A_{\theta,\eta}h_{\theta,\eta})A_{\theta,\eta}h_{\theta,\eta} = 0$ by the projection principle and is usually in the form of a conditional expectation. Section 2 in [22] provides an excellent guideline for searching for a least favorable submodel. Since the projection $\Pi_{\theta,\eta}\dot{\ell}_{\theta,\eta}$ on the closed linear span of the nuisance scores is not necessarily a nuisance score itself, the least favorable submodel may not always exist. However, we assume that in our setting a least favorable submodel always exists or can be approximated sufficiently closely by an *approximately least favorable submodel*. An insightful review of least favorable submodels and efficient score functions can be found in Chapter 3 of [13]. Systematic coverage of semiparametric efficiency theory can be found in [1] and [2].

The least favorable submodel in this paper will be constructed in the following manner. We consider a general map from the neighborhood of $\theta$ into the parameter set for $\eta$, denoted by $t \mapsto \eta_t(\theta, \eta)$. Then, the map $t \mapsto \ell(t, \theta, \eta)(x)$ can be defined as follows:

(1) $$\ell(t, \theta, \eta)(x) = \log \mathrm{lik}(t, \eta_t(\theta, \eta))(x).$$

The details of this map will depend on the situation.

2.2. *Notation and assumptions.* The dependence on $x \in \mathcal{X}$ of the likelihood and score quantities will be largely suppressed for clarity in this section and hereafter. The $\dot{\ell}(t,\theta,\eta)$, $\ddot{\ell}(t,\theta,\eta)$ and $\ell^{(3)}(t,\theta,\eta)$ are separately the first, second and third derivatives of $\ell(t,\theta,\eta)$ with respect to $t$. For brevity, we write $\dot{\ell}_0 = \dot{\ell}(\theta_0,\theta_0,\eta_0)$, $\ddot{\ell}_0 = \ddot{\ell}(\theta_0,\theta_0,\eta_0)$ and $\ell_0^{(3)} = \ell^{(3)}(\theta_0,\theta_0,\eta_0)$, where $\theta_0$ and $\eta_0$ are the true values of $\theta$ and $\eta$, respectively. Based on the definition of the least favorable submodel, $\dot{\ell}_0$ is just $\tilde{\ell}_0$ defined above. $\ell_\theta(t,\theta,\eta)$ indicates the first derivative of $\ell(t,\theta,\eta)$ w.r.t. $\theta$. Similarly, $\ell_{t,\theta}(t,\theta,\eta)$ denotes the derivative of $\dot{\ell}(t,\theta,\eta)$ w.r.t. $\theta$. Also, $\ell_{t,t}(\theta)$ and $\ell_{t,\theta}(\eta)$ indicate the maps $\theta \mapsto \ddot{\ell}(t,\theta,\eta)$ and $\eta \mapsto \ell_{t,\theta}(t,\theta,\eta)$, respectively. Let $\varrho_n$ denote $(\theta - \hat{\theta}_n)\tilde{I}_0^{1/2}$ and let $\phi(\cdot)$ ($\Phi(\cdot)$) represent the density (cumulative distribution) of a standard normal random variable. $\gtrsim$ and $\lesssim$ mean greater than, or smaller than, up to a universal constant. Define $x \vee y$ ($x \wedge y$) to be the maximum (minimum) value of $x$ and $y$.

$\mathbb{P}_n$ and $\mathbb{G}_n$ are used to denote the empirical distribution and the empirical process of the observations, respectively. Furthermore, we use the operator notation for evaluating expectation. Thus, for every measurable function $f$ and true probability $P$,

$$\mathbb{P}_n f = \frac{1}{n}\sum_{i=1}^n f(X_i), \qquad Pf = \int f\,dP \quad \text{and} \quad \mathbb{G}_n f = \frac{1}{\sqrt{n}}\sum_{i=1}^n (f(X_i) - Pf).$$

We now make the following assumptions:

1. $\theta_0 \in \Theta \subset \mathbb{R}^1$, where $\Theta$ is a compact set and $\theta_0$ is an interior point of $\Theta$;
2. $\eta_\theta(\theta,\eta) = \eta$ for any $(\theta,\eta) \in \Theta \times \mathcal{H}$;
3. $\|\hat{\eta}_{\tilde{\theta}_n} - \eta_0\| = O_P(n^{-1/2} + |\tilde{\theta}_n - \theta_0|)$ when $\tilde{\theta}_n = \theta_0 + o_P(1)$ for some norm $\|\cdot\|$;
4. the maps
$$(t,\theta,\eta) \mapsto \frac{\partial^{l+m}}{\partial t^l\, \partial \theta^m}\ell(t,\theta,\eta)$$
have integrable envelope functions in $L_1(P)$ in some neighborhood of $(\theta_0, \theta_0, \eta_0)$ for $(l,m) = (0,0), (1,0), (2,0), (3,0), (1,1), (1,2), (2,1)$;
5. there exists some neighborhood $V$ of $(\theta_0,\theta_0,\eta_0)$ in $\Theta \times \Theta \times \mathcal{H}$ such that the classes of functions $\{\ddot{\ell}(t,\theta,\eta)(x) : (t,\theta,\eta) \in V\}$ and $\{\ell_{t,\theta}(t,\theta,\eta)(x) : (t,\theta,\eta) \in V\}$ are $P$-Donsker and $\{\ell^{(3)}(t,\theta,\eta)(x) : (t,\theta,\eta) \in V\}$ is $P$-Glivenko–Cantelli;
6. 

(2) $$\mathbb{G}_n(\dot{\ell}(\theta_0,\theta_0,\eta) - \dot{\ell}_0) = O_P(\|\eta - \eta_0\|),$$

(3) $$P\ddot{\ell}(\theta_0,\theta_0,\eta) - P\ddot{\ell}(\theta_0,\theta_0,\eta_0) = O(\|\eta - \eta_0\|),$$

(4) $$P\ell_{t,\theta}(\theta_0,\theta_0,\eta) - P\ell_{t,\theta}(\theta_0,\theta_0,\eta_0) = O(\|\eta - \eta_0\|),$$

(5) $$P\dot{\ell}(\theta_0,\theta_0,\eta) = O(\|\eta - \eta_0\|^2)$$

for all $\eta$ in some neighborhood of $\eta_0$;
7. $\tilde{I}_0$ is strictly positive.

Assumption 2 ensures that the least favorable submodel passes through $(\theta, \eta)$, that is, $\ell(\theta, \theta, \eta)(x) = \log \operatorname{lik}(\theta, \eta)(x)$. Assumption 3 implicitly assumes that we have a metric or topology defined on the set of possible values of the nuisance parameter $\eta$. In this paper, uniform and weak topology norms are applied to the nuisance parameter in different examples. Definitions of the uniform and weak topology norms will be given in Section 5. Furthermore, the parametric convergence rate of the nuisance parameter is needed to obtain our second-order results. Assumption 4 can be viewed as comprising regular smoothness conditions on the Euclidean parameters of the least favorable submodel. Assumption 4 implies that $\ell(t, \theta, \eta)$ is smooth enough in its Euclidean parameter arguments so that $-P\dot{\ell}_0^2 = P\ddot{\ell}_0$. Assumption 4 also implies that $(\partial/\partial\theta)P\dot{\ell}(\theta_0, \theta, \eta_0) = 0$ at $\theta = \theta_0$. Fixing $\eta$ and differentiating $P_{\theta,\eta}\dot{\ell}(\theta, \theta, \eta)$ relative to $\theta$ gives $P_{\theta,\eta}\dot{\ell}_{\theta,\eta}\dot{\ell}(\theta, \theta, \eta) + P_{\theta,\eta}\ddot{\ell}(\theta, \theta, \eta) + (\partial/(\partial t))|_{t=\theta}P_{\theta,\eta}\dot{\ell}(\theta, t, \eta) = 0$ since $P_{\theta,\eta}\dot{\ell}(\theta, \theta, \eta) = 0$ for every $(\theta, \eta)$ and we can choose $(\theta, \eta) = (\theta_0, \eta_0)$.

The assumptions also impose some regular smoothness conditions on $\ell(t, \theta, \eta)$ relative to $\eta$ in the function space. Condition (2) involves the continuity modulus of the empirical processes. It can be easily satisfied if we can show that $\dot{\ell}(\theta_0, \theta_0, \eta) - \dot{\ell}_0$ divided by $\|\eta - \eta_0\|$ belongs to a $P$-Donsker class. The verification methods for (3)–(5) vary for different situations. Assuming a uniform norm is applied, (3) and (4) are usually satisfied if $\ell_{t,\theta}(\eta)$ and $\ell_{t,t}(\eta)$ have bounded Fréchet derivatives.

To verify (5), we need to briefly introduce Taylor series in Banach spaces [32]. Let $\zeta$ be a map from $\mathbb{D}_\zeta \subset \mathbb{D} \mapsto \mathbb{E}$, where $\mathbb{D}$ and $\mathbb{E}$ are both Banach spaces. If we assume $\zeta(\cdot)$ is second-order Fréchet differentiable, then the Taylor expansion of $\zeta(\vartheta + h)$ around $\zeta(\vartheta)$ can be written as $\zeta(\vartheta + h) = \zeta(\vartheta) + \zeta'_\vartheta(h) + \zeta''_{\vartheta+\tau h}(h, h)/2$, where $\tau \in [0, 1]$. $\zeta'_\vartheta(h)$ is just the regular Fréchet derivative of $\zeta(\cdot)$ at the point $\vartheta$ along the direction $h$ and $\zeta''_\vartheta(h, g)$ is the second-order Fréchet derivative from $\mathbb{D}_\zeta^2$ to $\mathbb{E}$. We can then write $P\dot{\ell}(\theta_0, \theta_0, \eta) = P[\frac{p_0 - p_{\theta_0,\eta}}{p_0}(\dot{\ell}(\theta_0, \theta_0, \eta) - \dot{\ell}(\theta_0, \theta_0, \eta_0))] - P[\dot{\ell}(\theta_0, \theta_0, \eta_0)(\frac{p_{\theta_0,\eta} - p_0}{p_0} - A_0(\eta - \eta_0))]$, where $A_0 = A_{\theta_0,\eta_0}$ and $A_{\theta,\eta}$ is the score operator for $\eta$ at $(\theta, \eta)$, for example, the Fréchet derivative of $\log p_{\theta,\eta}$ relative to $\eta$. The above equation holds since $P\tilde{\ell}_0 A_0 h = 0$ for every $h$, by the orthogonality property of the efficient score function. Note that the boundedness property of $\zeta''_\theta(\cdot, \cdot)$ means that $\|\zeta''_\theta(h, g)\|_\mathbb{E} \leq \|h\|_{\mathbb{D}_\zeta} \|g\|_{\mathbb{D}_\zeta}$. Thus, under the given regularity conditions, Fréchet differentiability of $\eta \mapsto \dot{\ell}(\theta_0, \theta_0, \eta)$ plus second-order Fréchet differentiability of $\eta \mapsto \operatorname{lik}(\theta_0, \eta)$ implies (5) based on the above discussions if the uniform norm is being applied to $\eta$.

In principle, these smoothness conditions on the least favorable submodel make the profile likelihood $pl_n(\theta)$ behave asymptotically like a parametric likelihood. The imposed assumptions are stronger than assumptions (3.1)–(3.4) in [22], enabling us to achieve higher order asymptotic expansions for the log-profile likelihood.

**3. Second-order asymptotic inference.** In this section, we present second-order asymptotic expansions of the log-profile likelihood which prepare us for deriving the main results of Section 4 on the higher order structure of the posterior profile distribution. Some of the results of this section are useful in their own right for inference about $\theta$. The assumptions of Section 2 are assumed throughout. We need the following lemma on the behavior of $\tilde{\theta}_n$, a random sequence of approximations of $\hat{\theta}_n$:

LEMMA 1. *If $\tilde{\theta}_n$ satisfies $(\tilde{\theta}_n - \hat{\theta}_n) = o_P(1)$, then*

$$(6) \qquad \mathbb{P}_n \dot{\ell}(\theta_0, \tilde{\theta}_n, \hat{\eta}_{\tilde{\theta}_n}) = \mathbb{P}_n \tilde{\ell}_0(X_i) + O_P(n^{-1/2} + |\tilde{\theta}_n - \hat{\theta}_n|)^2,$$

$$(7) \qquad \mathbb{P}_n \ddot{\ell}(\theta_0, \tilde{\theta}_n, \hat{\eta}_{\tilde{\theta}_n}) = P\ddot{\ell}_0 + O_P(n^{-1/2} + |\tilde{\theta}_n - \hat{\theta}_n|).$$

REMARK 1. Conditions (6) and (7) can essentially be viewed as the empirical versions of the no-bias conditions for the least favorable submodel (see, i.e., Chapter 25 of [28]). We can easily verify (6) and (7) if every argument of $\dot{\ell}(t, \theta, \eta)$ and $\ddot{\ell}(t, \theta, \eta)$ is smooth enough and the above empirical process assumptions are satisfied.

The following theorem gives key higher order expansions of the log-profile likelihood around $\hat{\theta}_n$ and on the error term in the asymptotic linearity expansion of $\hat{\theta}_n$.

THEOREM 1. *If $\tilde{\theta}_n$ satisfies $(\tilde{\theta}_n - \hat{\theta}_n) = o_P(1)$, then*

$$(8) \qquad \log pl_n(\tilde{\theta}_n) = \log pl_n(\hat{\theta}_n) - \frac{n}{2}(\tilde{\theta}_n - \hat{\theta}_n)^2 \tilde{I}_0 \\ + O_P(n|\tilde{\theta}_n - \hat{\theta}_n|^3 + n^{-1/2}),$$

$$(9) \qquad \sqrt{n}(\hat{\theta}_n - \theta_0) = \frac{1}{\sqrt{n}} \sum_{i=1}^{n} \tilde{\ell}_0(X_i) \tilde{I}_0^{-1} + O_P(n^{-1/2}).$$

REMARK 2. Expansions (8) and (9) are essentially second-order versions of (6) and (5), respectively, in [22], which have the respective error terms $o_P(\sqrt{n}|\tilde{\theta}_n - \theta_0| + 1)^2$ and $o_P(1)$. The parametric counterparts to (9) can be found in [16].

REMARK 3. Expansion (8) can be used to construct an estimator of the standard error of $\hat{\theta}_n$, which is called the "discretized" version of the observed profile information, $\hat{I}_n$, in [21]. Specifically, the discretized version of the observed profile information is expressed as a discretized second derivative of the profile likelihood in $\hat{\theta}_n$ as follows:

$$\hat{I}_n = -2\frac{\log pl_n(\hat{\theta}_n + s_n) - \log pl_n(\hat{\theta}_n)}{ns_n^2}. \tag{10}$$

Expansion (8) implies that

$$\hat{I}_n = \tilde{I}_0 + O_P(|s_n| + n^{-3/2}|s_n|^{-2}). \tag{11}$$

Obviously, the theoretically optimal step size of $\hat{I}_n$ is $s_n = O_P(n^{-1/2})$ and $s_n^{-1} = O_P(n^{1/2})$ in terms of the order of error term. In that case, $\hat{I}_n$ is a $\sqrt{n}$-consistent estimator of $\tilde{I}_0$.

An advantage of the method given in Remark 3 is that we can estimate $\tilde{I}_0$ even without an explicit form for the efficient Fisher information matrix or efficient score function. We only need the form of the profile likelihood, which is the minimal requirement, to carry out this numerical differentiation. Formula (11) provides us insight into the relationship between the step size of numerical differentiation and the convergence rate of $\hat{I}_n$. In other words, we can set a specific step size in advance to achieve the desired convergence rate. This is an improvement on Corollary 3 given in [21] which can only prove the consistency of the observed profile information.

**4. Main results.** We now present the main results on the posterior profile distribution. Let $\tilde{P}_{\theta|\tilde{X}}$ be the posterior profile distribution of $\theta$ w.r.t. the prior $\rho(\theta)$ given data $\tilde{X} = (X_1, \ldots, X_n)$. Define $\Delta_n(\theta) = n^{-1}\{\log pl_n(\theta) - \log pl_n(\hat{\theta}_n)\}$. A preliminary result, Theorem 2 with Corollaries 1 and 2 below, shows that the normal approximation to the posterior is second-order accurate for the cumulative distribution, the density and for the moments. The main result, Theorem 3, shows that the posterior profile distribution can be used to achieve second-order accurate frequentist inference.

THEOREM 2. *Assume the above assumptions and that*

$$\Delta_n(\tilde{\theta}_n) = o_P(1) \qquad \text{implies } \tilde{\theta}_n = \theta_0 + o_P(1). \tag{12}$$

*If proper prior $\rho(\theta_0) > 0$ and $\rho(\cdot)$ has a continuous and finite first-order derivative in some neighborhood of $\theta_0$, then we have, for $-\infty < \xi < \infty$,*

$$\sup_{\xi \in \mathbb{R}^1} |\tilde{P}_{\theta|\tilde{X}}(\sqrt{n}(\theta - \hat{\theta}_n)\tilde{I}_0^{1/2} \leq \xi) - \Phi(\xi)| = O_P(n^{-1/2}). \tag{13}$$

We note that the general theory concerning asymptotic expansions of posterior distributions in parametric models can be found in [11]. We also note that Theorem 1 in [17] implies the following:

$$\tilde{P}_{\theta|\tilde{X}}(\sqrt{n}(\theta - \hat{\theta}_n)\tilde{I}_0^{1/2} \leq \xi) = \Phi(\xi) + o_P(1). \tag{14}$$

Clearly, (14) is a first-order version of (13). A possibly more practical version of (13) is

$$\tilde{P}_{\theta|\tilde{X}}(\sqrt{n}(\theta - \hat{\theta}_n)\hat{I}_n^{1/2} \leq \xi) = \Phi(\xi) + O_P(n^{-1/2}), \tag{15}$$

where $\hat{I}_n$ can be estimated using (11) with an appropriate step size, for example, $s_n = O_P(n^{-1/2})$ and $s_n^{-1} = O_P(n^{1/2})$. Thus, we can construct the one-sided/two-sided credible set for $\theta$ with probability coverage $\alpha + O_P(n^{-1/2})$ in the following. Denoting by $z_\alpha$ to be the standard normal $\alpha$th quantile, we have

$$\tilde{P}_{\theta|\tilde{X}}\left(\theta \leq \hat{\theta}_n + \frac{z_\alpha}{\sqrt{nI}}\right) = \alpha + O_P(n^{-1/2}), \tag{16}$$

$$\tilde{P}_{\theta|\tilde{X}}\left(\hat{\theta}_n - \frac{z_{1-\alpha/2}}{\sqrt{nI}} \leq \theta \leq \hat{\theta}_n + \frac{z_{1-\alpha/2}}{\sqrt{nI}}\right) = 1 - \alpha + O_P(n^{-1/2}) \tag{17}$$

for $\alpha \in (0,1)$, where $I = \tilde{I}_0$ or $\hat{I}_n$.

COROLLARY 1. *Under the assumptions of Theorem 2, let $f_n(\cdot)$ be the posterior profile density of $\sqrt{n}\varrho_n$ relative to the prior $\rho(\theta)$. We then have*

$$f_n(\xi) = \phi(\xi) + O_P(n^{-1/2}). \tag{18}$$

REMARK 4. The parametric analog of (18) is (2.2) in [6], which is a higher order expansion of the multivariate posterior density of the vector $\sqrt{n}(\theta - \hat{\theta}_n)$ in a parametric model. Note that the parametric version involves the full likelihood rather than the profile likelihood and thus a prior is assigned to each element of the multivariate $\theta$. However, the posterior distributions relative to the full likelihood and the profile likelihood coincide for certain special priors which will be discussed in Remark 7 below.

COROLLARY 2. *Under the assumptions of Theorem 2 and recalling that $\varrho_n = (\theta - \hat{\theta}_n)\tilde{I}_0^{1/2}$, we have that if $\int_{-\infty}^{+\infty} |\theta|^r \rho(\theta)\,d\theta < \infty$, then*

$$\tilde{E}_{\theta|\tilde{X}}\varrho_n^r = n^{-r/2}EU^r + O_P(n^{-(r+1)/2}), \tag{19}$$

*where $\tilde{E}_{\theta|\tilde{X}}\varrho_n^r$ is the rth posterior moment of $\varrho_n$ and $U \sim N(0,1)$.*

REMARK 5. Note that the $r$th posterior moment of $\varrho_n$ in the above is based on the posterior profile distribution. By Corollary 2, we thus have

$$\hat{\theta}_n = \tilde{E}_{\theta|\tilde{X}}(\theta) + O_P(n^{-1}), \tag{20}$$

$$\tilde{I}_0 = \frac{1}{n\tilde{\text{Var}}_{\theta|\tilde{X}}(\theta)} + O_P(n^{-1/2}), \tag{21}$$

where

$$\tilde{\text{Var}}_{\theta|\tilde{X}}(\theta) = \tilde{E}_{\theta|\tilde{X}}(\theta - \tilde{E}_{\theta|\tilde{X}}(\theta))^2.$$

From (20), we know the maximum likelihood estimator of $\theta$ can be estimated by the mean of the profile sampler with an error of order $O_P(n^{-1})$. Moreover, from the proof in Section 7 of Theorem 3 below, we can verify that $\hat{\theta}_n$ is also estimated by the median of the profile sampler to the same order of accuracy. Similarly, the efficient information can be estimated by the inverse of the variance of the profile sampler with an error of order $O_P(n^{-1/2})$. This is a better method to estimate $\tilde{I}_0$ than (11) since it is automatically $\sqrt{n}$-consistent. Note that the first-order versions of (20) and (21) can be derived from Theorem 1 of [17].

Combining (9) and (20), we know that the mean value of the profile sampler can be shown to be a semiparametric efficient estimator of $\theta$. This conclusion also holds for the median value of the profile sampler. In this paper, we have provided an alternative efficient estimator to the maximum likelihood estimator $\hat{\theta}_n$.

We now present the main theorem of this paper. The $\alpha$th quantile of the posterior profile distribution, $\tau_{n\alpha}$, is defined as $\tau_{n\alpha} = \inf\{\xi : \tilde{P}_{\theta|\tilde{X}}(\theta \leq \xi) \geq \alpha\}$. Without loss of generality, $\tilde{P}_{\theta|\tilde{X}}(\theta \leq \tau_{n\alpha}) = \alpha$. We can also define $\kappa_{n\alpha} \equiv \sqrt{n}(\tau_{n\alpha} - \hat{\theta}_n)$, that is, $\tilde{P}_{\theta|\tilde{X}}(\sqrt{n}(\theta - \hat{\theta}_n) \leq \kappa_{n\alpha}) = \alpha$. The following theorem ensures that there exists a $\hat{\kappa}_{n\alpha}$ based on the data such that $P(\sqrt{n}(\hat{\theta}_n - \theta_0) \leq \hat{\kappa}_{n\alpha}) = \alpha$ and $|\hat{\kappa}_{n\alpha} - \kappa_{n\alpha}| = O_P(n^{-1/2})$.

THEOREM 3. *Under the assumptions of Theorem 2 and assuming that $\tilde{\ell}_0(X)$ has finite third moment with a nondegenerate distribution, there exists a $\hat{\kappa}_{n\alpha}$ based on the data such that $P(\sqrt{n}(\hat{\theta}_n - \theta_0) \leq \hat{\kappa}_{n\alpha}) = \alpha$ and $\hat{\kappa}_{n\alpha} - \kappa_{n\alpha} = O_P(n^{-1/2})$.*

REMARK 6. Note that the nondegenerate distribution assumption of $\tilde{\ell}(X)$ can be easily satisfied if $X$ has a nonsingular absolutely continuous component. Theorem 3 implies that the one- (two-) sided confidence interval for $\theta$ can be estimated by the one- (two-) sided credible set of the same level from the profile sampler with an error of the order $O_P(n^{-1})$. We conjecture that $\sqrt{n}$ times the $O_P(n^{-1/2})$ term in Theorem 3 converges to the product of two different nontrivial but uniformly integrable Gaussian processes.

REMARK 7. We can essentially generate the profile sampler from the marginal posterior of $\theta$ with respect to a certain joint prior on $\psi = (\theta, \eta)$ which is possibly data-dependent [17]. For example, in the Cox model with right-censored data, a gamma process prior on $\eta$ [12] with jumps at observed event times, but not involving $\theta$, can be such a prior.

**5. Examples.** We now illustrate the verification of the assumptions of Section 2 with three examples. The detailed technical illustrations and model assumptions for the three examples can be found in [19, 20, 21, 22, 23]. We also present simulation studies to assess the properties of the profile sampler for the first example.

5.1. *The Cox model with right-censored data.*

5.1.1. *Theory.* The Cox model is

$$(22) \quad \lambda(t|z) \equiv \lim_{\Delta \to 0} \frac{1}{\Delta} Pr(t \leq T < t + \Delta | T \geq t, Z = z) = \lambda(t) \exp(\theta z),$$

where $\lambda$ is an unspecified baseline hazard function and $\theta$ is a vector including the regression parameters [4]. For the Cox model applied to right-censored failure time data, we observe that $X = (Y, \delta, Z)$, where $Y = T \wedge C$, $\delta = I\{T \leq C\}$ and $Z \in \mathbb{Z} \subset \mathbb{R}^1$ is a regression covariate. $T$ is a failure time with integrated hazard $e^{\theta z}\Lambda(t)$ given the covariate $Z$, where $\Lambda(y) = \int_0^y \lambda(t)\,dt$ is a cadlag, monotone increasing cumulative hazard function with $\Lambda(0) = 0$. $C$ is a censoring time independent of $T$ given $Z$. We define a likelihood for the parameter $(\theta, \Lambda)$ by replacing $\lambda(y)$ with the point mass $\Lambda\{y\}$:

$$(23) \quad \text{lik}(\theta, \Lambda) = (e^{\theta z}\Lambda\{y\}e^{-e^{\theta z}\Lambda(y)})^\delta (e^{-e^{\theta z}\Lambda(y)})^{1-\delta}.$$

$\theta$ is assumed to come from some compact set $\Theta$ and the true regression coefficient, $\theta_0$, belongs to the interior of $\Theta$. The parameter space for $\Lambda$, $\mathcal{H}$, is restricted to a set of nondecreasing, cadlag functions on the interval $[0, \tau]$, with $\Lambda(\tau) \leq M$ for a given constant $M$.

We now discuss the form of the profile likelihood. Suppose there are $l$ observed failures at times $T_{(1)} < \cdots < T_{(l)}$, where $(i)$ is the label for the $i$th ordered failure and $t_i$ is the observed value of $T_{(i)}$. $z_{[i]}$ is the covariate corresponding to the observed event time $t_i$. The log-profile likelihood (equivalently, the log-partial likelihood) for $\theta$ is given by

$$(24) \quad \log pl_n(\theta) = \sum_{i=1}^{l}\left(\theta z_{[i]} - \log \sum_{j \in R_i} e^{\theta z_j}\right),$$

where $R_i = \{j : Y_j \geq t_i\}$ is the risk set. In this case, the profiled nuisance parameter is not present in $pl_n(\theta)$. Nevertheless, it is not hard to verify that

$$\hat{\Lambda}_\theta(t) = \sum_{\{Y_i \leq t\}} \frac{\delta_i}{\sum_{j \in R_i} \exp(\theta z_j)}. \tag{25}$$

Note that $\hat{\Lambda}_\theta$ is a nondecreasing step function with support points at the observed event times and, based on [10], $\|\hat{\Lambda}_{\tilde{\theta}_n} - \Lambda_0\|_\infty = O_P(n^{-1/2} + |\tilde{\theta}_n - \theta_0|)$.

The score function for $\theta$ can be easily derived as

$$\dot{\ell}_{\theta,\Lambda}(x) = \delta z - z e^{\theta z} \Lambda(y).$$

Given a fixed $\Lambda$ and a bounded function $h : \mathbb{R}^1 \mapsto \mathbb{R}^1$, we can define a path $\Lambda_t$ by $d\Lambda_t(y) = (1 + th(y)) d\Lambda(y)$. Thus, the score function for $\Lambda$ in the direction $h$ via an operator $A_{\theta,\Lambda} : L_2(\Lambda) \mapsto L_2(P_{\theta,\Lambda})$ is $A_{\theta,\Lambda} h(y, \delta, z) = \delta h(y) - e^{\theta z} \int_{[0,y]} h \, d\Lambda$. Following the regular conditions and discussions on page 16 of [22], the least favorable direction $h_{\theta,\Lambda}$ at $(\theta, \Lambda)$ can be constructed as

$$h_{\theta,\Lambda}(y) = \frac{E_{\theta,\Lambda} e^{\theta Z} Z 1\{Y \geq y\}}{E_{\theta,\Lambda} e^{\theta Z} 1\{Y \geq y\}}.$$

Substituting $\theta = t$ and $\Lambda = \Lambda_t(\theta, \Lambda)$ [where $d\Lambda_t(\theta, \Lambda) = (1 + (\theta - t)h_0) d\Lambda$ and $h_0(\cdot)$ is an abbreviation for $h_{\theta_0,\Lambda_0}(\cdot)$] in the above Cox likelihood and differentiating with respect to $t$, we obtain

$$\dot{\ell}(t, \theta, \Lambda)(x) = \dot{\ell}_{t,\Lambda_t(\theta,\Lambda)} - A_{t,\Lambda_t(\theta,\Lambda)}\left(\frac{h_0(y)}{1 + (\theta - t)h_0(y)}\right)(x),$$

$$= \delta z - z e^{tz} \Lambda_t(\theta, \Lambda)(y) - \delta \frac{h_0(y)}{1 + (\theta - t)h_0(y)} + e^{tz} \int_0^y h_0 \, d\Lambda,$$

$$\ddot{\ell}(t, \theta, \Lambda)(x) = -\delta \frac{h_0^2(y)}{(1 + (\theta - t)h_0(y))^2}$$

$$- z^2 e^{tz} \Lambda_t(\theta, \Lambda)(y) + 2z e^{tz} \int_0^y h_0 \, d\Lambda,$$

$$\ell_{t,\theta}(t, \theta, \Lambda)(x) = -z e^{tz} \int_0^y h_0 \, d\Lambda + \delta \frac{h_0(y)^2}{(1 + (\theta - t)h_0(y))^2},$$

$$\ell^{(3)}(t, \theta, \Lambda)(x) = -2\delta \frac{h_0^3(y)}{(1 + (\theta - t)h_0(y))^3}$$

$$- z^3 e^{tz} \Lambda_t(\theta, \Lambda)(y) + 3z^2 e^{tz} \int_0^y h_0 \, d\Lambda.$$

We know the maps $(t, \theta, \Lambda) \mapsto \ell^{(k)}(t, \theta, \Lambda)$, for $k = 1, 2, 3$, are continuous and uniformly bounded around $(\theta_0, \theta_0, \Lambda_0)$, relative to the uniform topology on $\Lambda$, by the inequality $\int_0^y h_0 \, d(\Lambda - \Lambda_0) \lesssim \|\Lambda - \Lambda_0\|_\infty \|h_0\|_{BV}$, where $\|h_0\|_{BV}$ is the total variation of $h_0(\cdot)$ in $[0, \tau]$. The total variation of a function $f : [a, b] \mapsto R$, $\|f\|_{BV}$, is $|f(a)| + \int_{(a,b]} |df(s)|$. Considering the fact that the class of uniformly bounded functions with bounded variation over compacta is $P$-Donsker, we can check the empirical process assumptions by repeatedly using the $P$-Donsker preservation results. The following lemmas verify the remaining conditions, thus the results of Sections 3 and 4 hold.

LEMMA 2. *Under the above set-up for the proportional hazards Cox model, assumption* 6 *is satisfied.*

LEMMA 3. *Under the above set-up for the proportional hazards Cox model, condition* (12) *is satisfied.*

5.1.2. *Simulation study.* To verify that the profile sampler can generate second-order frequentist valid inference, we conducted simulations for Cox regression with right-censored data for various sample sizes under a Lebesgue prior. For each sample size, 500 data sets were analyzed. The event times were generated from (22) with one covariate $Z \sim U[0, 1]$. The regression coefficient is $\theta = 1$ and $\Lambda(t) = \exp(t) - 1$. The censoring time $C \sim U[0, t_n]$, where $t_n$ was chosen such that the average effective sample size over 500 samples is approximately $0.9n$. For each dataset, Markov chains of length 5,000 with a burn-in period of 1,000 were generated using the Metropolis algorithm. The jumping density for the coefficient was normal with current iteration and variance tuned to yield an acceptance rate of 20–40%. The approximate variance of the estimator of $\theta$ was computed by numerical differentiation with step size proportional to $n^{-1/2}$, according to Remark 3.

Table 1 summarizes the results from the simulations giving the average across 500 samples of the maximum likelihood estimate (MLE), mean of the profile sampler (CM), mean squared difference between two estimates of $\theta$ ($\text{MSD}_E$), estimated standard errors based on MCMC ($\text{SE}_M$), estimated standard errors based on numerical derivatives ($\text{SE}_N$), mean squared difference between the two estimated standard errors ($\text{MSD}_V$) and empirical coverage of nominal 0.95 confidence intervals based on MCMC (CP95). The Monte Carlo standard error of CP95 is $\approx 0.01 = \sqrt{0.05 \times 0.95/500}$. Table 2 summarizes the difference of boundaries for the two-sided 95% confidence interval for $\theta$ generated by numerical differentiation, that is, (17), and MCMC, respectively. $\text{LB}_M$ ($\text{LB}_N$) and $\text{UB}_M$ ($\text{UB}_N$) denote the lower and upper bound, respectively, of the confidence interval from the MCMC chain (numerical derivative).

In all cases, the bias in Table 1 is small. Similar simulations are also performed in the Cox model with current status data in [17], that is, Table 1, which has larger bias. The contrast of two simulations reveals an interesting

Table 1
*Simulation results for Cox regression with right-censored data based on 500 samples (the true value of the regression coefficient is $\theta = 1$)*

| $n$ | MLE | CM | $\sqrt{\mathrm{MSD_E}}$ | $\mathrm{SE_M}$ | $\mathrm{SE_N}$ | $\sqrt{\mathrm{MSD_V}}$ | CP95 |
|---|---|---|---|---|---|---|---|
| 20  | 1.1049 | 1.1376 | 0.0978 | 4.4128 | 4.3004 | 0.2934 | 0.9496 |
| 50  | 1.0202 | 1.0262 | 0.0275 | 3.9869 | 3.9548 | 0.1378 | 0.9496 |
| 100 | 1.0156 | 1.0181 | 0.0195 | 3.8592 | 3.8561 | 0.1568 | 0.9506 |
| 200 | 1.0131 | 1.0147 | 0.0114 | 3.8124 | 3.8105 | 0.1220 | 0.9500 |
| 500 | 1.0012 | 1.0016 | 0.0069 | 3.7598 | 3.7691 | 0.1206 | 0.9502 |

$n$, sample size; MLE, maximum likelihood estimator; CM, empirical mean; $\mathrm{MSD_E}$, mean squared difference between two estimates of $\theta$; $\mathrm{SE_M}$, estimated standard errors based on MCMC; $\mathrm{SE_N}$, estimated standard errors based on numerical derivatives; $\mathrm{MSD_V}$, mean squared difference between the two estimated standard errors; CP95, empirical coverage of nominal 0.95 confidence intervals based on MCMC.

phenomenon: the profile sampler based on the semiparametric models with faster convergence rate is more accurate. Note that the terms $n|\mathrm{MLE} - \mathrm{CM}|$, $\sqrt{n}|\mathrm{SE_M} - \mathrm{SE_N}|$, $n|\mathrm{LB_M} - \mathrm{LB_N}|$ and $n|\mathrm{UB_M} - \mathrm{UB_N}|$ are bounded in probability according to Corollary 1 and Theorem 3, that is, (11), (20) and (21). The realizations of these terms summarized in Table 2 clearly illustrate their boundedness. Based on the above results, we can conclude that the profile sampler is a second-order frequentist valid procedure.

5.2. *The proportional odds model with right-censored data.* The survival function in this example is parameterized such that the ratios of the odds of survival for subjects with different covariates are constant with time: the conditional survival function $S_Z(u)$ of the event time, $T$, given the covariate $Z$ satisfies $-\mathrm{logit}(S_Z(u)) = \log \eta(u) + Z\theta$, where $\mathrm{logit}(y) = \log(y/(1-y))$.

Table 2
*Simulation results for confidence intervals*

| $n$ | $n|\mathrm{MLE} - \mathrm{CM}|$ | $\sqrt{n}|\mathrm{SE_M} - \mathrm{SE_N}|$ | $n|\mathrm{LB_M} - \mathrm{LB_N}|$ | $n|\mathrm{UB_M} - \mathrm{UB_N}|$ |
|---|---|---|---|---|
| 20  | 0.6541 | 0.5027 | 0.1920 | 2.2823 |
| 50  | 0.3062 | 0.2270 | 0.1809 | 1.1212 |
| 100 | 0.2587 | 0.0311 | 0.5987 | 0.1301 |
| 200 | 0.3218 | 0.0279 | 0.4810 | 0.5253 |
| 500 | 0.2017 | 0.2080 | 0.7524 | 0.3518 |

$\mathrm{LB_M}$ ($\mathrm{UB_M}$), lower (upper) bound of the 95% confidence interval based on MCMC; $\mathrm{LB_N}$ ($\mathrm{UB_N}$), lower (upper) bound of the 95% confidence interval based on numerical derivative.

We define the likelihood as

$$\text{lik}(\theta, \eta) = \left[\frac{e^{-z\theta}\eta\{y\}}{(\eta(y) + e^{-z\theta})(\eta(y-) + e^{-z\theta})}\right]^\delta \left[\frac{e^{-z\theta}}{\eta(y) + e^{-z\theta}}\right]^{1-\delta}, \quad (26)$$

where $\eta\{y\}$ is the jump size in $\eta$ at $y$. The score function for $\theta$ is

$$\dot{\ell}_{\theta,\eta}(x) = -z\left(1 - \frac{e^{-z\theta}}{\eta(y) + e^{-z\theta}} - \frac{\delta e^{-z\theta}}{\eta(y-) + e^{-z\theta}}\right).$$

The score function for $\eta$ via the direction of bounded function $h \in L_2(\eta)$ is

$$A_{\theta,\eta}h(x) = (\partial/\partial t)|_{t=0}\ell_{\theta,\eta_t} = \delta h(y) - \frac{\int_0^y h\, d\eta}{\eta(y) + e^{-z\theta}} - \frac{\delta \int_0^{y-} h\, d\eta}{\eta(y-) + e^{-z\theta}},$$

where $d\eta_t = (1+th)\, d\eta$. Let $A^*_{\theta,\eta}$ denote the adjoint of $A_{\theta,\eta}$. Then, $A^*_{\theta,\eta}A_{\theta,\eta}h(u)$ is the information operator for the nuisance parameter $\eta$ when $\theta$ is known. It is shown to be continuously invertible on the space of functions of bounded variation on $[0,\tau]$ in Lemma 4.3 of [19]. Hence, the least favorable direction $h_0$ is defined as $(A^*_{\theta_0,\eta_0}A_{\theta_0,\eta_0})^{-1}A^*_{\theta_0,\eta_0}\ell_{\theta_0,\eta_0}$. The form of the information operator and $A^*_{\theta,\eta}$ can be found in [19]. By setting $d\eta_t(\theta, \eta) = (1+(\theta-t)h_0)\, d\eta$, we can obtain $\ell(t, \theta, \eta) = \log \text{lik}(t, \eta_t(\theta, \eta))$. Hence, the maps in assumption 4 can be derived as follows:

$$\dot{\ell}(t, \theta, \eta)(x) = -z\left(1 - \frac{e^{-zt}}{\eta_t(y) + e^{-zt}} - \frac{\delta e^{-zt}}{\eta_t(y-) + e^{-zt}}\right)$$
$$- \delta \frac{h_0(y)}{1 + (\theta-t)h_0(y)} + \frac{\int_0^y h_0\, d\eta}{\eta_t(y) + e^{-zt}} + \frac{\delta \int_0^{y-} h_0\, d\eta}{\eta_t(y-) + e^{-zt}},$$

$$\ddot{\ell}(t, \theta, \eta)(x) = -\delta \frac{h_0(y)^2}{(1 + (\theta-t)h_0(y))^2} + W^a_{t,\theta,\eta}(y, z) + \delta W^a_{t,\theta,\eta}(y-, z),$$

$$W^a_{t,\theta,\eta}(y, z) = \frac{(\int_0^y h_0\, d\eta + ze^{-zt})^2 - z^2 e^{-zt}(\eta_t(y) + e^{-zt})}{(\eta_t(y) + e^{-zt})^2},$$

$$\ell_{t,\theta}(t, \theta, \eta)(x) = \delta \frac{h_0^2(y)}{(1 + (\theta-t)h_0(y))^2} + W^b_{t,\theta,\eta}(y, z) + \delta W^b_{t,\theta,\eta}(y-, z),$$

$$W^b_{t,\theta,\eta}(y, z) = -\frac{\int_0^y h_0\, d\eta(ze^{-zt} + \int_0^y h_0\, d\eta)}{(\eta_t(y) + e^{-zt})^2},$$

$$\ell^{(3)}(t, \theta, \eta) = -2\frac{\delta h_0^3(y)}{(1 + (\theta-t)h_0(y))^3} + W^c_{t,\theta,\eta}(y, z) + \delta W^c_{t,\theta,\eta}(y-, z),$$

$$W^c_{t,\theta,\eta}(y, z) = \frac{2(ze^{-zt} + \int_0^y h_0\, d\eta)^3}{(\eta_t(y) + e^{-zt})^3}$$
$$- \frac{z^2 e^{-zt}(e^{-zt} + \eta_t(y))(2ze^{-zt} + 3\int_0^y h_0\, d\eta - z\eta_t(y))}{(\eta_t(y) + e^{-zt})^3}.$$

Under the regular conditions, the above maps are continuous and uniformly bounded around $(\theta_0, \theta_0, \eta_0)$ by the same reasoning as was used in the first example. We also know that $\|\hat{\eta}_{\tilde{\theta}_n} - \eta_0\|_\infty = O_P(n^{-1/2} + |\tilde{\theta}_n - \theta_0|)$ by Theorem 3.1 in [21]. By similar techniques to those used in the first example, we can easily verify assumption 5. The following lemmas verify the remaining conditions.

LEMMA 4. *Under the above set-up for the proportional odds model, assumption 6 is satisfied.*

LEMMA 5. *Under the above set-up for the proportional odds model, condition (12) is satisfied.*

5.3. *Case-control studies with a missing covariate.* The third example is a logistic regression model for case-control studies with a missing covariate considered by [23] and [24]. We observe two independent random samples of sizes $n_C$ and $n_R$ from the distributions of $(D, W, Z)$ and $(D, W)$, respectively. Following the assumptions concerning the distribution of the random vector $(D, W, Z)$ in [24], we can construct the likelihood for the vector $(D, W, Z)$ in the form $p_\theta(d, w|z) \, d\eta(z)$, where

$$(27) \quad p_\theta(d, w|z) = (\Xi_{\gamma,\beta}(z))^d (1 - \Xi_{\gamma,\beta}(z))^{1-d} \frac{1}{\sigma} \phi\left(\frac{\omega - \alpha_0 - \alpha_1 z}{\sigma}\right),$$

$\Xi_{\gamma,\beta}(z) = (1 + \exp(-\gamma - \beta e^z))^{-1}$ and $d\eta$ denotes the density of $\eta$ with respect to some dominating measure on $\mathcal{Z} \subset \mathbb{R}^1$.

We assume $n_C = n_R$ so that the observations can then be paired. Here, we denote the complete sample components by $Y_C = (D_C, W_C)$ and $Z_C$ and the reduced sample components by $Y_R = (D_R, W_R)$. Thus, the likelihood is defined as

$$(28) \quad \text{lik}(\theta, \eta)(x) = p_\theta(y_C|z_C)\eta\{z_C\} \int p_\theta(y_R|z) \, d\eta(z).$$

The unknown parameters are $\theta = (\beta, \alpha_0, \alpha_1, \gamma, \sigma)$, ranging over a compact $\Theta \subset \mathbb{R}^4 \times (0, \infty)$, and the distribution $\eta$ of the regression variable which is restricted to the set of nondegenerate probability distributions with support within a known compact interval. We will concentrate on the regression coefficient $\beta$, considering $\theta_2 = (\alpha_0, \alpha_1, \gamma, \sigma)$ and $\eta$ as nuisance parameters.

We start by introducing the least favorable submodel. The score function of $\theta$, $\dot{\ell}_{\theta,\eta}$, is the summation of the score functions for the conditional density $p_\theta(y_C|z_C)$ and that for the mixture density $p_\theta(y_R|\eta)$, given as follows:

$$\dot{\ell}_\theta(y_C|z_C) = \frac{\partial}{\partial \theta} \log p_\theta(y_C|z_C) \quad \text{and} \quad \dot{\ell}_{\theta,\eta}(y_R) = \frac{\int \dot{\ell}_\theta(y_R|z) p_\theta(y_R|z) \, d\eta(z)}{p_\theta(y_R|\eta)}.$$

Furthermore, by defining $d\eta_t = (1 + th)\,d\eta$, where $h$ is an arbitrary bounded function satisfying $\int h\,d\eta = 0$, we can obtain the score function for $\eta$ in the direction $h$,

$$A_{\theta,\eta}h(x) = h(z_C) + \frac{\int h(z) p_\theta(y_R|z)\,d\eta(z)}{p_\theta(y_R|\eta)}.$$

By the projection principle discussed in [24], we thus define the least favorable submodel as follows:

$$\ell(t, \beta, \theta_2, \eta) = \log l(\theta_t(\theta, \eta), \eta_t(\theta, \eta)),$$

where $\theta_t(\theta, \eta) = \theta - (\beta - t)a_0$, $d\eta_t(\theta, \eta) = (1 + (\beta - t)a_0^T(h_0 - \eta h_0))\,d\eta$ and $a_0^T = (1, -\tilde{I}_{0,12}(\tilde{I}_{0,22})^{-1})$. The efficient information matrix $\tilde{I}_0$ can be decomposed into four submatrices corresponding to parameters $\beta$ and the group $(\alpha_0, \alpha_1, \gamma, \sigma)$. $\tilde{I}_{0,ij}$ corresponds to the $(i,j)$th block of $\tilde{I}_0$ for $i = 1, 2$ and $j = 2$. In Section 8 of [23], the least favorable direction at the true value, $h_0$, is proved to be a bounded and Lipschitz continuous function.

Let $(\hat{\theta}_{2,\beta}, \hat{\eta}_\beta)$ be the profile likelihood estimator for $(\theta_2, \eta)$ when $\beta$ is given so that $\hat{\theta}_\beta = (\beta, \hat{\theta}_{2,\beta})$. [21] showed that

$$(29) \qquad \|\hat{\eta}_{\tilde{\beta}_n} - \eta_0\|_{BL_1} + \|\hat{\theta}_{\tilde{\beta}_n} - \theta_0\| = O_P(|\tilde{\beta}_n - \beta_0| + n^{-1/2})$$

for any $\tilde{\beta}_n$ consistent for $\beta_0$. The norm applied to the function $\eta$ and vector $\theta$ is the weak topology norm and Euclidean norm, respectively. The weak topology norm on $\eta$ is defined as $\|\eta\|_{BL_1} = \sup_{h \in BL_1} |\int h(z)\,d\eta(z)|$, where $BL_1$ denotes the set of all functions $h : \mathcal{Z} \mapsto [-1, 1]$ that are Lipschitz norm bounded above by 1, that is, $|h(z_1) - h(z_2)| \leq \|z_1 - z_2\|_{\mathcal{Z}}$. The following lemmas verify the remaining conditions.

LEMMA 6. *Under the above set-up for the case-control model, assumptions* 4–6 *are satisfied.*

LEMMA 7. *Under the above set-up for the case-control model, condition* (12) *is satisfied.*

**6. Discussion.** Our theory ensures second-order frequentist correctness of the profile Bayes analysis for the finite-dimensional parameter. The necessary and sufficient conditions required for third or higher order frequentist inference need to be constructed in order to complete general higher order semiparametric frequentist inference theory in the future. Our future work could also include extending our methods to semiparametric models with slower convergence rates for the nuisance parameter, for example, $\|\hat{\eta}_{\tilde{\theta}_n} - \eta_0\| = O_P(n^{-1/3} + \|\tilde{\theta}_n - \theta_0\|)$, as happens with the Cox model for current status data. The conjecture in Remark 6 implies that the $O_P(n^{-1/2})$ rate in Theorem 2 is sharp. Hence, to show this conjecture may be a future research goal, although it appears to be very challenging.

## 7. Proofs.

PROOF OF LEMMA 1. We first show the following no-bias conditions:

$$P\dot{\ell}(\theta_0, \tilde{\theta}_n, \hat{\eta}_{\tilde{\theta}_n}) = O_P(n^{-1/2} + |\tilde{\theta}_n - \theta_0|)^2, \tag{30}$$

$$P\ddot{\ell}(\theta_0, \tilde{\theta}_n, \hat{\eta}_{\tilde{\theta}_n}) = P\ddot{\ell}_0 + O_P(n^{-1/2} + |\tilde{\theta}_n - \theta_0|). \tag{31}$$

(30) can be written as $P\dot{\ell}(\theta_0, \tilde{\theta}_n, \hat{\eta}_{\tilde{\theta}_n}) - P\dot{\ell}_0 = [P\dot{\ell}(\theta_0, \tilde{\theta}_n, \hat{\eta}_{\tilde{\theta}_n}) - P\dot{\ell}(\theta_0, \theta_0, \hat{\eta}_{\tilde{\theta}_n})] + [P\dot{\ell}(\theta_0, \theta_0, \hat{\eta}_{\tilde{\theta}_n}) - P\dot{\ell}_0]$. The second square bracket is bounded by $O_P(\|\hat{\eta}_{\tilde{\theta}_n} - \eta_0\|^2)$, by (5). By the ordinary two-term Taylor expansion, the first square bracket equals

$$(\tilde{\theta}_n - \theta_0)(\partial/\partial\theta)|_{\theta=\theta_0} P\dot{\ell}(\theta_0, \theta, \hat{\eta}_{\tilde{\theta}_n})$$
$$+ (1/2)(\tilde{\theta}_n - \theta_0)^2 \times (\partial^2/\partial\theta^2)|_{\theta=\theta^*} P\dot{\ell}(\theta_0, \theta, \hat{\eta}_{\tilde{\theta}_n}),$$

where $\theta^*$ is an intermediate value between $\tilde{\theta}_n$ and $\theta_0$. The second term of this expansion is of order $|\tilde{\theta}_n - \theta_0|^2$, by assumption 4. We now consider the first term. Define $\eta \mapsto L(\eta) = (\partial/\partial\theta)|_{\theta=\theta_0} P\dot{\ell}(\theta_0, \theta, \eta)$. Then, $L(\hat{\eta}_{\tilde{\theta}_n}) - L(\eta_0) = O_P(\|\hat{\eta}_{\tilde{\theta}_n} - \eta_0\|)$, by (4) in assumption 6. Combining this with the fact that $L(\eta_0) = 0$, we have $(\tilde{\theta}_n - \theta_0) \times (\partial/\partial\theta)|_{\theta=\theta_0} P\dot{\ell}(\theta_0, \theta, \hat{\eta}_{\tilde{\theta}_n}) = O_P(n^{-1/2} + |\tilde{\theta}_n - \theta_0|)^2$. This completes the proof of (30). By assumption 3, the smoothness conditions on $\ddot{\ell}(t, \theta, \eta)$ and (3), we can also show (31) using similar analysis.

Recall that $\dot{\ell}_0(X) = \tilde{\ell}_0(X)$. It then suffices to show (6) if $\mathbb{G}_n\sqrt{n}(\dot{\ell}(\theta_0, \tilde{\theta}_n, \hat{\eta}_{\tilde{\theta}_n}) - \dot{\ell}_0) = O_P(\sqrt{n}|\tilde{\theta}_n - \theta_0| + 1)$. Note that, by (2), $\mathbb{G}_n\sqrt{n}(\dot{\ell}(\theta_0, \tilde{\theta}_n, \hat{\eta}_{\tilde{\theta}_n}) - \dot{\ell}_0) = \sqrt{n}(\tilde{\theta}_n - \theta_0)\mathbb{G}_n\ell_{t,\theta}(\theta_0, \theta_n^*, \hat{\eta}_{\tilde{\theta}_n}) + \sqrt{n}O_P(\|\hat{\eta}_{\tilde{\theta}_n} - \eta_0\|)$, where $\theta_n^*$ is an intermediate value between $\theta_0$ and $\tilde{\theta}_n$. Combining this with assumption 5, we have proven (6). Considering assumption 5 and (31), we can prove (7). □

PROOF OF THEOREM 1. We first show (9). Note that $0 = \mathbb{P}_n\dot{\ell}(\hat{\theta}_n, \hat{\theta}_n, \hat{\eta}_n) = \mathbb{P}_n\dot{\ell}(\theta_0, \hat{\theta}_n, \hat{\eta}_n) + (\hat{\theta}_n - \theta_0)\mathbb{P}_n\ddot{\ell}(\theta_0, \hat{\theta}_n, \hat{\eta}_n) + ((\hat{\theta}_n - \theta_0)^2/2)\mathbb{P}_n\ell^{(3)}(\theta_n^*, \hat{\theta}_n, \hat{\eta}_n)$, where $\theta_n^*$ is intermediate between $\theta_0$ and $\hat{\theta}_n$. By considering Lemma 1 and assumption 5, we construct the following equation about $(\hat{\theta}_n - \theta_0)$: $0 = n^{-1}\sum_{i=1}^n \tilde{\ell}_0(x_i) + (\hat{\theta}_n - \theta_0)P\ddot{\ell}_0 + O_P(n^{-1})$. This completes the proof of (9). To prove (8), we first show that

$$\log pl_n(\tilde{\theta}_n) = \log pl_n(\theta_0) + (\tilde{\theta}_n - \theta_0)\sum_{i=1}^n \tilde{\ell}_0(X_i) - \frac{n}{2}(\tilde{\theta}_n - \theta_0)^2 \tilde{I}_0 \tag{32}$$
$$+ O_P(n|\tilde{\theta}_n - \hat{\theta}_n|^3 + n^{-1/2})$$

for any $\tilde{\theta}_n$ satisfying $(\tilde{\theta}_n - \hat{\theta}_n) = o_P(1)$. Note that

$$n^{-1}(\log pl_n(\tilde{\theta}_n) - \log pl_n(\theta_0)) = \mathbb{P}_n\ell(\tilde{\theta}_n, \tilde{\theta}_n, \hat{\eta}_{\tilde{\theta}_n}) - \mathbb{P}_n\ell(\theta_0, \theta_0, \hat{\eta}_{\theta_0}).$$

The right-hand side of the above equation is bounded below and above by $\mathbb{P}_n(\ell(\tilde{\theta}_n, \tilde{\psi}_n) - \ell(\theta_0, \tilde{\psi}_n))$, where the lower and upper bound separately correspond to $\tilde{\psi}_n = (\theta_0, \hat{\eta}_{\theta_0})$ and $(\tilde{\theta}_n, \hat{\eta}_{\tilde{\theta}_n})$. We then apply a three-term Taylor expansion to both upper and lower bounds. By considering Lemma 1 and assumption 5, we find that the upper bound and the lower bound match at the order of $O_P(n^{-3/2} + |\tilde{\theta}_n - \hat{\theta}_n|^3)$. We have thus proven (32). By replacing $\tilde{\theta}_n$ with $\hat{\theta}_n$ in (32), we have

$$\log pl_n(\hat{\theta}_n) = \log pl_n(\theta_0) + (\hat{\theta}_n - \theta_0) \sum_{i=1}^n \tilde{\ell}_0(X_i) - \frac{n}{2}(\hat{\theta}_n - \theta_0)^2 \tilde{I}_0$$
(33)
$$+ O_P(n^{-1/2}).$$

The difference between (32) and (33) gives (8) by considering (9). □

PROOF OF THEOREM 2. Suppose that $F_n(\cdot)$ is the posterior profile distribution of $\sqrt{n}\varrho_n$ w.r.t. the prior $\rho(\theta)$, where $\varrho_n = (\theta - \hat{\theta}_n)\tilde{I}_0^{1/2}$. The whole proof of Theorem 2 can be briefly summarized in the following expression:

$$F_n(\xi) = \frac{\int_{-\infty}^{\xi n^{-1/2}} \rho(\hat{\theta}_n + \varrho_n \tilde{I}_0^{-1/2})(pl_n(\hat{\theta}_n + \varrho_n \tilde{I}_0^{-1/2}))/(pl_n(\hat{\theta}_n))\, d\varrho_n}{\int_{-\infty}^{+\infty} \rho(\hat{\theta}_n + \varrho_n \tilde{I}_0^{-1/2})(pl_n(\hat{\theta}_n + \varrho_n \tilde{I}_0^{-1/2}))/(pl_n(\hat{\theta}_n))\, d\varrho_n}.$$

For the denominator, we first prove that the posterior mass outside $|\varrho_n| \leq r_n$ is of arbitrarily small order, where $r_n = o(n^{-1/3})$ and $\sqrt{n}r_n \to \infty$. The mass inside this integration region can be approximated by a stochastic polynomial in powers of $n^{-1/2}$ with an error of the order $O_P(n^{-1})$. The numerator can be analyzed similarly. Finally, the asymptotic expansions of both numerator and denominator yield the quotient series, which is the desired result. We first state some lemmas before the giving formal proof of Theorem 2.

LEMMA 2.1. *Let $r_n = o(n^{-1/3})$ and $\sqrt{n}r_n \to \infty$. Under the conditions of Theorem 2, we have*

(34) $$\int_{|\varrho_n|>r_n} \rho(\hat{\theta}_n + \varrho_n \tilde{I}_0^{-1/2}) \frac{pl_n(\hat{\theta}_n + \varrho_n \tilde{I}_0^{-1/2})}{pl_n(\hat{\theta}_n)}\, d\varrho_n = O_P(n^{-1}).$$

PROOF. Fix $r > 0$. We then have

$$\int_{|\varrho_n|>r} \rho(\hat{\theta}_n + \varrho_n \tilde{I}_0^{-1/2}) \frac{pl_n(\hat{\theta}_n + \varrho_n \tilde{I}_0^{-1/2})}{pl_n(\hat{\theta}_n)}\, d\varrho_n$$
$$\leq I\{\Delta_n^r < -n^{-1/2}\} \exp(-\sqrt{n}) \int_\Theta \rho(\theta)\, d\theta + I\{\Delta_n^r \geq -n^{-1/2}\},$$

20    G. CHENG AND M. R. KOSOROK

where $\Delta_n^r = \sup_{|\varrho_n|>r} \Delta_n(\hat{\theta}_n + \varrho_n \tilde{I}_0^{-1/2})$. By a minor revision of Lemma A.1 in the Appendix of [17], we have $I\{\Delta_n^r > -n^{-1/2}\} = O_P(n^{-1})$. This implies that there exists a positive decreasing sequence $r_n = o(n^{-1/3})$ with $\sqrt{n}r_n \to \infty$ such that (34) holds. □

LEMMA 2.2.  *Let $r_n = o(n^{-1/3})$ and $\sqrt{n}r_n \to \infty$. Under the conditions of Theorem 2, we have*

$$(35) \quad \int_{-r_n}^{r_n} \left| \frac{pl_n(\hat{\theta}_n + \varrho_n \tilde{I}_0^{-1/2})}{pl_n(\hat{\theta})} \rho(\hat{\theta}_n + \varrho_n \tilde{I}_0^{-1/2}) - \exp\left(-\frac{n}{2}\varrho_n^2\right) \rho(\hat{\theta}_n) \right| d\varrho_n$$
$$= O_P(n^{-1}).$$

PROOF.  The posterior mass over the region $|\varrho_n| \leq r_n$ is bounded by

$$(*) \quad \int_{|\varrho_n| \leq r_n} \left| \frac{pl_n(\hat{\theta}_n + \varrho_n \tilde{I}_0^{-1/2})}{pl_n(\hat{\theta})} \rho(\hat{\theta}_n) - \exp\left(-\frac{n}{2}\varrho_n^2\right) \rho(\hat{\theta}_n) \right| d\varrho_n$$

$$(**) \quad + \int_{|\varrho_n| \leq r_n} \left| \frac{pl_n(\hat{\theta}_n + \varrho_n \tilde{I}_0^{-1/2})}{pl_n(\hat{\theta})} \rho(\hat{\theta}_n + \varrho_n \tilde{I}_0^{-1/2}) \right.$$
$$\left. - \frac{pl_n(\hat{\theta}_n + \varrho_n \tilde{I}_0^{-1/2})}{pl_n(\hat{\theta})} \rho(\hat{\theta}_n) \right| d\varrho_n.$$

Using (8), we obtain

$$(*) = \int_{|\varrho_n| \leq r_n} \left[ \rho(\hat{\theta}_n) \exp\left(-\frac{n\varrho_n^2}{2}\right) |\exp(O_P(n|\varrho_n|^3 + n^{-1/2})) - 1| \right] d\varrho_n$$

$$= n^{-1/2} \int_{|u_n| \leq \sqrt{n}r_n} \left[ \rho(\hat{\theta}_n) \exp\left(-\frac{u_n^2}{2}\right) \right.$$
$$\left. \times |\exp(n^{-1/2}(|u_n|^3 + 1)O_P(1)) - 1| \right] du_n$$

$$= n^{-1} \times O_P(1) \times \int_{|u_n| \leq \sqrt{n}r_n} \left[ \rho(\hat{\theta}_n) \exp\left(-\frac{u_n^2}{2}\right) (|u_n|^3 + 1) \right] du_n$$

$$= O_P(n^{-1}),$$

where the second equality follows by replacing $\sqrt{n}\varrho_n$ with $u_n$ and the third equality follows from the fact that $|\exp(n^{-1/2}(|u_n|^3 + 1)O_P(1)) - 1| = O_P(1)n^{-1/2}(|u_n|^3 + 1)$ since $|u_n| \leq \sqrt{n}r_n$ and $r_n = o(n^{-1/3})$, that is, $u_n = o(n^{1/6})$. By the following analysis of $(**)$, we can also show $(**) = O_P(n^{-1})$

since $\exp(O_P(n\varrho_n^3 + n^{-1/2})) = O_P(1)$ with $|\varrho_n| \leq r_n$:

$$(**) = \int_{|\varrho_n| \leq r_n} \left[ |\varrho_n \tilde{I}_0^{-1/2} \dot{\rho}(\theta_n^*)| \exp\left(-\frac{n}{2}\varrho_n^2 + O_P(n\varrho_n^3 + n^{-1/2})\right) \right] d\varrho_n$$

$$\leq M \int_{|\varrho_n| \leq r_n} \left[ |\varrho_n| \exp\left(-\frac{n}{2}\varrho_n^2\right) \right] d\varrho_n \times \sup_{|\varrho_n| \leq r_n} \exp(O_P(n\varrho_n^3 + n^{-1/2})),$$

where $\theta_n^*$ is an intermediate value between $\hat{\theta}_n$ and $\hat{\theta}_n + \varrho_n \tilde{I}_0^{-1/2}$. □

We next start the formal proof of Theorem 2. First, note that

$$\int_{-\infty}^{+\infty} \left[ \rho(\hat{\theta}_n + \varrho_n \tilde{I}_0^{-1/2}) \frac{pl_n(\hat{\theta}_n + \varrho_n \tilde{I}_0^{-1/2})}{pl_n(\hat{\theta}_n)} \right] d\varrho_n$$

$$= \int_{|\varrho_n| \geq r_n} \left[ \rho(\hat{\theta}_n + \varrho_n \tilde{I}_0^{-1/2}) \frac{pl_n(\hat{\theta}_n + \varrho_n \tilde{I}_0^{-1/2})}{pl_n(\hat{\theta}_n)} \right] d\varrho_n$$

$$+ \int_{|\varrho_n| \leq r_n} \left[ \rho(\hat{\theta}_n + \varrho_n \tilde{I}_0^{-1/2}) \frac{pl_n(\hat{\theta}_n + \varrho_n \tilde{I}_0^{-1/2})}{pl_n(\hat{\theta}_n)} \right] d\varrho_n.$$

By Lemma 2.1, the first integral on the right-hand side is of the order $O_P(n^{-1})$. The second integral on the right-hand side can be decomposed into the following summands:

$$\int_{|\varrho_n| \leq r_n} \left[ \rho(\hat{\theta}_n + \varrho_n \tilde{I}_0^{-1/2}) \frac{pl_n(\hat{\theta}_n + \varrho_n \tilde{I}_0^{-1/2})}{pl_n(\hat{\theta})} - \exp\left(-\frac{n}{2}\varrho_n^2\right) \rho(\hat{\theta}_n) \right] d\varrho_n$$

$$+ \int_{|\varrho_n| \leq r_n} \left[ \exp\left(-\frac{n}{2}\varrho_n^2\right) \rho(\hat{\theta}_n) \right] d\varrho_n.$$

The first part is bounded by $O_P(n^{-1})$ via Lemma 2.2. The second part equals

$$n^{-1/2} \rho(\hat{\theta}_n) \int_{|u_n| \leq \sqrt{n} r_n} e^{-u_n^2/2} \, du_n = n^{-1/2} \rho(\hat{\theta}_n) \int_{-\infty}^{+\infty} e^{-u_n^2/2} \, du_n + O(n^{-1}),$$

where $u_n = \sqrt{n}\varrho_n$. The above equality follows from the inequality that $\int_x^\infty e^{-y^2/2} \, dy \leq x^{-1} e^{-x^2/2}$ for any $x > 0$.

Consolidating the above analysis, we have

(36)
$$\int_{-\infty}^{+\infty} \left[ \rho(\hat{\theta}_n + \varrho_n \tilde{I}_0^{-1/2}) \frac{pl_n(\hat{\theta}_n + \varrho_n \tilde{I}_0^{-1/2})}{pl_n(\hat{\theta}_n)} \right] d\varrho_n$$
$$= n^{-1/2} \rho(\hat{\theta}_n) \sqrt{2\pi} + O_P(n^{-1})$$

and, by similar analysis, we obtain

$$
(37) \quad \int_{-\infty}^{\xi n^{-1/2}} \left[ \rho(\hat{\theta}_n + \varrho_n \tilde{I}_0^{-1/2}) \frac{pl_n(\hat{\theta}_n + \varrho_n \tilde{I}_0^{-1/2})}{pl_n(\hat{\theta}_n)} \right] d\varrho_n
$$
$$
= n^{-1/2} \rho(\hat{\theta}_n) \int_{-\infty}^{\xi} e^{-y^2/2} \, dy + O_P(n^{-1}).
$$

The quotient of (36) and (37) generates the desired result, (13). This completes the proof of Theorem 2 in its entirety. □

PROOF OF COROLLARY 1. From the proof of Theorem 2, we have

$$
\tilde{P}_{\theta|\tilde{X}}(\sqrt{n}(\theta - \hat{\theta}_n)\tilde{I}_0^{1/2} \leq \xi)
$$
$$
= \frac{\int_{-\infty}^{\xi n^{-1/2}} \rho(\hat{\theta}_n + \varrho_n \tilde{I}_0^{-1/2})(pl_n(\hat{\theta}_n + \varrho_n \tilde{I}_0^{-1/2}))/(pl_n(\hat{\theta}_n)) \, d\varrho_n}{\int_{-\infty}^{+\infty} \rho(\hat{\theta}_n + \varrho_n \tilde{I}_0^{-1/2})(pl_n(\hat{\theta}_n + \varrho_n \tilde{I}_0^{-1/2}))/(pl_n(\hat{\theta}_n)) \, d\varrho_n}.
$$

By differentiating both sides relative to $\xi$ and combining with (36), we obtain

$$
f_n(\xi) = \frac{\rho(\hat{\theta}_n + \xi \tilde{I}_0^{-1/2}/\sqrt{n})(pl_n(\hat{\theta}_n + \xi \tilde{I}_0^{-1/2} n^{-1/2}))/(pl_n(\hat{\theta}_n))}{\sqrt{2\pi} \rho(\hat{\theta}_n) + O_P(n^{-1/2})}.
$$

Based on (8), the numerator in the above equals $\rho(\hat{\theta}_n) \exp(-\xi^2/2) + O_P(n^{-1/2})$ by some analysis. This completes the proof. □

PROOF OF COROLLARY 2. The expansion in (19) is the quotient of two expansions of the form (36) and (37). We can see this as follows. First,

$$
\tilde{E}_{\theta|x}(\varrho_n^r) = \frac{\int_{-\infty}^{+\infty} \varrho_n^r \rho(\hat{\theta}_n + \varrho_n \tilde{I}_0^{-1/2})(pl_n(\hat{\theta}_n + \varrho_n \tilde{I}_0^{-1/2}))/(pl_n(\hat{\theta}_n)) \, d\varrho_n}{\int_{-\infty}^{+\infty} \rho(\hat{\theta}_n + \varrho_n \tilde{I}_0^{-1/2})(pl_n(\hat{\theta}_n + \varrho_n \tilde{I}_0^{-1/2}))/(pl_n(\hat{\theta}_n)) \, d\varrho_n}.
$$

The denominator is $n^{-1/2}\sqrt{2\pi}\rho(\hat{\theta}_n) + O_P(n^{-1})$ by (36). Similarly, by the proof of Theorem 2, we know the numerator is $n^{-(r+1)/2}\rho(\hat{\theta}_n)\sqrt{2\pi}EU^r + O_P(n^{-(r+2)/2})$, equivalently, $(2/n)^{(r+1)/2}\Gamma((r+1)/2)\rho(\hat{\theta}_n) + O_P(n^{-(r+2)/2})$, where $U \sim N(0,1)$. Obviously, the quotient is $n^{-r/2}EU^r + O_P(n^{-(r+1)/2})$. If $r$ is odd, the quotient is simply $O_P(n^{-(r+1)/2})$. □

PROOF OF THEOREM 3. We first show that for any $\xi \in (0, \frac{1}{2})$ and $\xi < \alpha < 1 - \xi$,

$$
(38) \quad \tau_{n\alpha} = \hat{\theta}_n + \frac{z_\alpha}{\sqrt{n\tilde{I}_0}} + O_P(n^{-1}).
$$

Implicit in (13) is an expansion of $\tau_{n\alpha}$ in terms of $z_\alpha$. First, we set $\tau_{n\alpha} = \hat\theta_n + z_\alpha/\sqrt{n\tilde I_0} + r_n$ and we can then show $r_n = O_P(n^{-1})$. Plugging $\tau_{n\alpha}$ into (13), we obtain $\alpha = \tilde P_{\theta|\tilde X}(\sqrt{n}(\theta-\hat\theta_n)\tilde I_0^{1/2} \leq z_\alpha + \sqrt{n}\tilde I_0^{1/2} r_n) = \alpha + O_P(n^{-1/2}) + \sqrt{n}\tilde I_0^{1/2} r_n f_n(\tau_{n\alpha}^*)$, where $\tau_{n\alpha}^*$ is between $z_\alpha$ and $z_\alpha + \sqrt{n}\tilde I_0^{1/2} r_n$. The first equality comes from the definition of $\tau_{n\alpha}$. The second equality follows from Taylor expansion and (18). We can now deduce from these two equalities that $r_n = -O_P(n^{-1})\tilde I_0^{-1/2}(\Phi(\kappa_{n\alpha}^*) + O_P(n^{-1/2}))^{-1} = O_P(n^{-1})$ based on (18). Note that $r_n$ is well defined since $f_n(\tau_{n\alpha}^*)$ is strictly positive when $\xi < \alpha < 1-\xi$. This completes the proof of (38). Next the classical Edgeworth expansion implies that $P(n^{-1/2}\sum_{i=1}^n \tilde\ell_0(X_i)\tilde I_0^{-1/2} \leq z_\alpha + a_n) = \alpha$, where $a_n = O(n^{-1/2})$, for $\xi < \alpha < 1-\xi$. Let $\hat\kappa_{n\alpha} = z_\alpha \tilde I_0^{-1/2} + (\sqrt{n}(\hat\theta_n - \theta_0) - \frac{1}{\sqrt{n}}\sum_{i=1}^n \tilde\ell_0(X_i)\tilde I_0^{-1}) + a_n \tilde I_0^{-1/2}$. Then, $P(\sqrt{n}(\hat\theta_n - \theta_0) \leq \hat\kappa_{n\alpha}) = P(n^{-1/2}\sum_{i=1}^n \tilde\ell_0(X_i)\tilde I_0^{-1/2} \leq z_\alpha + a_n) = \alpha$. Combining (38) and (9), we obtain $\hat\kappa_{n\alpha} = \kappa_{n\alpha} + O_P(n^{-1/2})$. □

PROOF OF LEMMA 2. We first compute the Fréchet derivatives of $\ell_{t,\theta}(\Lambda)$ around $(\theta_0, \theta_0, \Lambda_0)$ by means of $\Lambda_s(y) = \Lambda(y) + s\int_0^y h\,d\Lambda = \Lambda(y) + sW_\Lambda(y)$, where $h(\cdot)$ is an arbitrary bounded function. The corresponding Fréchet derivatives are as follows:

$$\dot\ell_{t,\theta,\Lambda}(W_\Lambda) = -ze^{tz}\int_0^y h_0\,dW_\Lambda.$$

The operator $\dot\ell_{t,\theta,\Lambda}(W_\Lambda)$ is linear and continuous by the inequality

$$(39) \quad |\dot\ell_{t,\theta,\Lambda}(W_\Lambda) - \dot\ell_{t,\theta,\Lambda}(V_\Lambda)| \lesssim \left|\int_0^y h_0\,d(W_\Lambda - V_\Lambda)\right| \lesssim \|W_\Lambda - V_\Lambda\|_\infty,$$

almost surely, since $\Lambda$ is a cumulative hazard function with support $[0,\tau]$. It is also a bounded operator since we can replace $V_\Lambda$ with zero in (39). Note that $\dot\ell_{t,\theta,\Lambda}(0) = 0$ by its linearity. By similar reasoning, we can also know that $\ddot\ell_{t,t,\Lambda}(W_\Lambda)$ and $\ddot\ell_{t,\Lambda}(W_\Lambda)$ are both linear, continuous and bounded operators when $(t,\theta)$ is in some neighborhood of $(\theta_0, \theta_0)$ and $\Lambda \in \mathcal{H}$. The boundedness of the above two operators ensure that $P(\dot\ell_{t,\theta}(\theta_0,\theta_0,\Lambda) - \dot\ell_{t,\theta}(\theta_0,\theta_0,\Lambda_0)) = O_P(\|\Lambda - \Lambda_0\|_\infty)$ and $P(\ddot\ell(\theta_0,\theta_0,\Lambda) - \ddot\ell(\theta_0,\theta_0,\Lambda_0)) = O_P(\|\Lambda - \Lambda_0\|_\infty)$ when $\Lambda$ is in some neighborhood of $\Lambda_0$. To verify (5), we need to show that $\Lambda \mapsto \mathrm{lik}(\theta_0, \Lambda)$ is second-order Fréchet differentiable around $\Lambda_0$. To this end, the first derivative is

$$\dot{\mathrm{lik}}_\Lambda(W_\Lambda) = \frac{e^{\theta_0 z}}{\exp(e^{\theta_0 z}\Lambda(y))}\left(-\int_0^y dW_\Lambda\right)^{1-\delta}\left(W_\Lambda\{y\} - \Lambda\{y\}e^{\theta_0 z}\int_0^y dW_\Lambda\right)^\delta,$$

while the second derivative is

$$\ddot{\mathrm{lik}}_\Lambda(W_\Lambda, V_\Lambda) = \frac{e^{\theta_0 z}}{\exp(e^{\theta_0 z}\Lambda(y))}(e^{\theta_0 z}W_\Lambda(y)V_\Lambda(y))^{1-\delta}$$

$$\times (e^{2\theta_0 z}\Lambda\{y\}V_\Lambda(y)W_\Lambda(y)$$
$$- e^{\theta_0 z}W_\Lambda\{y\}V_\Lambda(y) - e^{\theta_0 z}V_\Lambda\{y\}W_\Lambda(y))^\delta.$$

Clearly, $\ddot{\text{lik}}_\Lambda(W_\Lambda, V_\Lambda)$ is a bounded bilinear operator. Its continuity follows from the continuity of the maps $W_\Lambda \mapsto \ddot{\text{lik}}_\Lambda(W_\Lambda, \cdot)$ and $V_\Lambda \mapsto \ddot{\text{lik}}_\Lambda(\cdot, V_\Lambda)$.

Next, we need to show that $\mathbb{G}_n(\dot{\ell}(\theta_0, \theta_0, \Lambda) - \dot{\ell}_0) = O_P(\|\Lambda - \Lambda_0\|_\infty)$. First, note that the class of functions of $(z, y)$,

$$\{ze^{\theta_0 z}(M\mathbf{1}\{\Lambda(y) \leq t\} - \Lambda_0(y)): \|\Lambda - \Lambda_0\|_\infty \leq \gamma, t \in [0, \tau]\},$$

is a VC class for each $\gamma < \infty$. Since $\Lambda$ is monotone and bounded by $M < \infty$, we now have that the class $\{ze^{\theta_0 z}(\Lambda(y) - \Lambda_0(y)): \|\Lambda - \Lambda_0\|_\infty \leq \gamma\}$ is a VC-hull class for each $\gamma < \infty$. Since $k_0 \gamma$ is an envelope for this last class for some $k_0 < \infty$ that does not depend on $\gamma$, we can use Theorem 2.14.1 in [29] to obtain that $\mathbb{G}_n z e^{\theta_0 z}(\Lambda(y) - \Lambda_0(y)) = O_P(\|\Lambda - \Lambda_0\|_\infty)$. A similar argument can be used to verify that $\mathbb{G}_n e^{\theta_0 z} \int_0^y h_0(s)(d\Lambda(s) - d\Lambda_0(s)) = O_P(\|\Lambda - \Lambda_0\|_\infty)$. Thus, $\mathbb{G}_n(\dot{\ell}(\theta_0, \theta_0, \Lambda) - \dot{\ell}_0) = O_P(\|\Lambda - \Lambda_0\|_\infty)$, as desired. □

PROOF OF LEMMA 3. The proof of Lemma 3 is analogous to that of Lemma 4 in [17], which is for the more general odds-rate model. □

PROOF OF LEMMA 4. The proof of Lemma 4 is analogous to that of Lemma 2. We can similarly verify the linearity, continuity and boundedness of $\ell_{t,\eta}(W_\eta)$, $\ell_{t,\theta,\eta}(W_\eta)$ and $\ell_{t,t,\eta}(W_\eta)$, whose concrete forms can be found in [3]. The verification of (5) also follows similar reasoning as used in the proof of Lemma 2. The forms of $\dot{\text{lik}}_\eta(W_\eta)$ and $\ddot{\text{lik}}_\eta(W_\eta, V_\eta)$ are specified in [3]. By analysis similar to that in the proof of Lemma 2, we can show (2). This completes the proof. □

PROOF OF LEMMA 5. The proof of Lemma 5 is analogous to that of Lemma 4 in [17], which is for the more general odds-rate model. □

PROOF OF LEMMA 6. Before we start the proof of Lemma 6, we first present the following necessary computations according to (28).

$$(40) \quad \dot{\ell}_\theta(y|z) = \begin{pmatrix} \frac{e^z}{\exp(\gamma+\beta e^z)+1} + (d-1)e^z \\ \frac{w-\alpha_0-\alpha_1 z}{\sigma^2} \\ \frac{z(w-\alpha_0-\alpha_1 z)}{\sigma^2} \\ \frac{1}{\exp(\gamma+\beta e^z)+1} + (d-1) \\ -\frac{1}{\sigma} + \frac{(w-\alpha_0-\alpha_1 z)^2}{\sigma^3} \end{pmatrix},$$

$$(41) \quad \ddot{\ell}_\theta(y|z) = \begin{pmatrix} -\frac{\exp(\gamma+\beta e^z)e^{2z}}{(1+\exp(\gamma+\beta e^z))^2} & 0 & 0 & -\frac{\exp(\gamma+\beta e^z)e^z}{(1+\exp(\gamma+\beta e^z))^2} & 0 \\ 0 & -\frac{1}{\sigma^2} & -\frac{z}{\sigma^2} & 0 & -\frac{2(w-\alpha_0-\alpha_1 z)}{\sigma^3} \\ 0 & 0 & -\frac{z^2}{\sigma^2} & 0 & -\frac{2z(w-\alpha_0-\alpha_1 z)}{\sigma^3} \\ 0 & 0 & 0 & -\frac{\exp(\gamma+\beta e^z)}{(1+\exp(\gamma+\beta e^z))^2} & 0 \\ 0 & 0 & 0 & 0 & \frac{1}{\sigma^2} - \frac{3(w-\alpha_0-\alpha_1 z)^2}{\sigma^4} \end{pmatrix},$$

where $\ddot{\ell}_\theta(y|z) = \partial^2/\partial\theta^2 \log p_\theta(y|z)$. We now compute $(t,\theta,\eta) \mapsto \partial^{l+m}/\partial t^l \partial\theta^m \ell(t,\theta,\eta)$, with the abbreviations $\theta_t = \theta_t(\theta,\eta)$ and $\eta_t = \eta_t(\theta,\eta)$, for $(l,m) = (0,0), (1,0), (2,0), (3,0), (1,1), (1,2), (2,1)$, as follows:

$$\dot{\ell}(t,\theta,\eta) = a_0^T(\dot{\ell}_{\theta_t}(y_C|z_C) + \dot{\ell}_{\theta_t,\eta_t}(y_R) - A_{\theta_t,\eta_t}(G_\eta^{h_0}(\beta,t))(x));$$

$$\ddot{\ell}(t,\theta,\eta) = a_0^T\left(\ddot{\ell}_{\theta_t}(y_C|z_C)a_0 + \frac{\partial}{\partial t}\dot{\ell}_{\theta_t,\eta_t}(y_R) - \frac{\partial}{\partial t}A_{\theta_t,\eta_t}(G_\eta^{h_0}(\beta,t))(x)\right);$$

$$\ell_{t,\beta}(t,\theta,\eta) = a_0^T\bigg(\ddot{\ell}_{\theta_t}(y_C|z_C)(\mathbf{1}_1 - a_0)$$
$$+ \frac{\partial}{\partial\beta}\dot{\ell}_{\theta_t,\eta_t}(y_R) - \frac{\partial}{\partial\beta}A_{\theta_t,\eta_t}(G_\eta^{h_0}(\beta,t))(x)\bigg);$$

$$\ell_{t,\theta_j}(t,\theta,\eta) = a_0^T\left(\ddot{\ell}_{\theta_t}(y_C|z_C)\mathbf{1}_j + \frac{\partial}{\partial\theta_j}\dot{\ell}_{\theta_t,\eta_t}(y_R) - \frac{\partial}{\partial\theta_j}A_{\theta_t,\eta_t}(G_\eta^{h_0}(\beta,t))(x)\right);$$

$$\ell^{(3)}(t,\theta,\eta) = \sum_i\sum_j\sum_k a_{ijk}\frac{\partial^3}{\partial\theta_k\partial\theta_j\partial\theta_i}\bigg|_{\theta=\theta_t}\log p_\theta + \sum_j\sum_i a_{ij}\frac{\partial}{\partial t}\ddot{\ell}_{\theta_t,\eta_t}\{ij\}$$
$$- \sum_j\sum_i a_{ij}\frac{\partial}{\partial t}\dot{A}_{\theta_t,\eta_t}\{ij\}.$$

For brevity, we omit the complete versions of the above formulas and refer the interested reader to [3]. However, we present the complete description of $\ddot{\ell}(t,\theta,\eta)$ in the following to illustrate some functional properties of the above formulas, which will be used in the proof of Lemma 6:

$$a_0^T\frac{\partial}{\partial t}\dot{\ell}_{\theta_t,\eta_t}(y_R) = a_0^T\bigg\{\frac{\eta_t(\ddot{\ell}_{\theta_t}p_{\theta_t}) + \eta_t(\dot{\ell}_{\theta_t}\dot{\ell}_{\theta_t}^T p_{\theta_t}) - \eta(\dot{\ell}_{\theta_t}H_0(z)^T p_{\theta_t})}{\eta_t p_{\theta_t}}$$
$$+ \frac{\eta_t(\dot{\ell}_{\theta_t}p_{\theta_t})}{\eta_t p_{\theta_t}} \times \frac{\eta(H_0(z)^T p_{\theta_t}) - \eta_t(\dot{\ell}_{\theta_t}^T p_{\theta_t})}{\eta_t p_{\theta_t}}\bigg\}a_0$$

and

$$a_0^T\frac{\partial}{\partial t}A_{\theta_t,\eta_t}(G_\eta^{h_0}(\beta,t))(x)$$
$$= a_0^T\bigg\{\frac{H_0(z)H_0(z)^T}{(1+(\beta-t)a_0^T H_0(z))^2} + \frac{\eta(H_0(z)\dot{\ell}_{\theta_t}^T p_{\theta_t})}{\eta_t p_{\theta_t}}$$
$$- \frac{\eta(H_0(z)p_{\theta_t})}{\eta_t p_{\theta_t}} \times \frac{\eta_t(\dot{\ell}^T p_{\theta_t}) - \eta(H_0(z)^T p_{\theta_t})}{\eta_t p_{\theta_t}}\bigg\}a_0.$$

In the preceding, $\mathbf{1}_i$ is a five-dimensional vector with the $i$th element one, and the others zero, $\theta_j$ is the $j$th element of the vector $\theta$, $a_i$ is the $i$th element of vector $a_0$, $a_{ij} \equiv a_i a_j$, $a_{ijk} \equiv a_i a_j a_k$ and $p_{\theta_t}$, $\dot{\ell}_{\theta_t}$ and $\ddot{\ell}_{\theta_t}$ are respective abbreviations for $p_{\theta_t}(y|z)$, $\dot{\ell}_{\theta_t}(y|z)$ and $\ddot{\ell}_{\theta_t}(y|z)$. $L^T$ is the transpose of $L$ and $H_0(z) \equiv h_0(z) - \eta h_0$. Note that $H_0(z) \in \mathcal{C}_1^1(\mathcal{Z})$ with zero mean after proper rescaling. $A_{\theta_t,\eta_t}(\cdot)$ is an abbreviation for $A_{\theta_t,\eta_t}((h_0 - \eta h_0)/(1+(\beta-t)a_0^T(h_0 - \eta h_0)))(x)$. $\ddot{\ell}_{\theta_t,\eta_t}\{ij\}$ and $\dot{A}_{\theta_t,\eta_t}\{ij\}$ are the respective $(i,j)$th elements of square matrices $(\partial/\partial t)\dot{\ell}_{\theta_t,\eta_t}(y_R)a_0^T(a_0a_0^T)^{-1}$ and $(\partial/\partial t)A_{\theta_t,\eta_t}(\cdot)a_0^T(a_0a_0^T)^{-1}$. The above notation is valid for $i,j,k = 1,\ldots,5$. We need the following two lemmas to verify assumption 4. □

LEMMA 6.1. *Given $z$ in some compact set $\mathcal{Z}$, $\theta$ and $\eta$ in some neighborhood of $\theta_0$ and $\eta_0$, respectively, we have*

$$\text{(42)} \qquad \frac{p_\theta(y|z)}{c_\theta(w) \int p_\theta(y|z)\,d\eta(z)} \in \mathcal{C}_1^1(\mathcal{Z}),$$

*where $c_\theta(w) = M_0 \exp(\frac{M|w|}{2\sigma^2})(|w|+1)$ and $0 < M_0, M < \infty$.*

PROOF. Note that $||w| - |\alpha_0 + \alpha_1 z|| \leq |w - \alpha_0 - \alpha_1 z| \leq |w| + |\alpha_0 + \alpha_1 z|$, and $z$ is in some compact set. Thus, we have the following inequalities:

$$\text{(43)} \qquad \exp\left(-\frac{|w|M_2}{2\sigma^2}\right) \lesssim \frac{p_\theta(y|z)}{\int p_\theta(y|z)\,d\eta(z)} \lesssim \exp\left(\frac{|w|M_1}{2\sigma^2}\right)$$

and

$$\text{(44)} \qquad \left|\frac{\partial}{\partial z}\left(\frac{p_\theta(y|z)}{\int p_\theta(y|z)\,d\eta(z)}\right)\right| \lesssim \exp\left(\frac{|w|M_1}{2\sigma^2}\right)(|w|+1),$$

where $M_i$ is some positive finite number, $i = 1, 2$. □

LEMMA 6.2. *Let $h_\theta(y|z) = \sum_{l=0}^L g_l(z;\sigma,\gamma,\beta)(w - \alpha_0 - \alpha_1 z)^l$ for $\theta \in \Theta$, where $g_l(z;\sigma,\gamma,\beta) \in \mathcal{C}_1^1(\mathcal{Z})$ and is continuous w.r.t. $\theta$ for $l = 0, 1, \ldots, L$. The following then has an integrable envelope function in $L_K(P)$ and is continuous at $(\theta, \eta_1, \eta_2)$ when $\theta$ is in some neighborhood of $\theta_0$ and $\eta_i$ is in some neighborhood of $\eta_0$ for $i = 1, 2$, and where $K$ is any positive integer:*

$$\text{(45)} \qquad f^h_{\theta,\eta_1,\eta_2}(y) \equiv \frac{\int h_\theta(y|z)p_\theta(y|z)\,d\eta_1(z)}{\int p_\theta(y|z)\,d\eta_2(z)}.$$

PROOF. The following is the envelope function for $f^h_{\theta,\eta_1,\eta_2}(y)$, $F^h(y)$:

$$|f^h_{\theta,\eta_1,\eta_2}(y)| \lesssim \sum_{l=0}^L (|w|+1)^l \frac{\int p_\theta(y|z)\,d\eta_1(z)}{\int p_\theta(y|z)\,d\eta_2(z)}$$

$$\lesssim \sum_{l=0}^L (|w|+1)^l \exp\left(\frac{|w|M_1}{2\sigma_{\min}^2}\right) \equiv F^h(y).$$

In the above, the first inequality follows from (45), the second one follows from (43) and $0 < \sigma_{\min} \leq \sigma \leq \sigma_{\max} < \infty$. Next, we only need to show $P|F^h(y)|^K < \infty$ for any positive integer $K$. Accordingly,

$$P|F^h_\theta(y)|^K \leq \sum_{i=0}^{1} \int_{\mathcal{Z}} \int_{-\infty}^{+\infty} \left(\sum_{l=0}^{L}(|w|+1)^l\right)^K$$
$$\times \exp\left(\frac{KM_1}{2\sigma_{\min}^2}|w|\right) p_{\theta_0}(w, d=i|z) \, dw \, d\eta_0(z)$$
$$\lesssim \sum_{i=0}^{1} \int_{\mathcal{Z}} \int_{-\infty}^{+\infty} \left(\sum_{l=0}^{L}(|w|+1)^l\right)^K \exp\left(\frac{KM_1}{2\sigma_{\min}^2}|w|\right)$$
$$\times \exp\left(-\frac{(|w| - |\alpha_0 + \alpha_1 z|)^2}{2\sigma_{\max}^2}\right) dw \, d\eta_0(z)$$
$$\lesssim \sum_{i=0}^{1} \int_{\mathcal{Z}} \int_{-\infty}^{+\infty} \left(\sum_{l=0}^{L}(|w|+1)^l\right)^K \exp\left(-\frac{(|w| - M_3)^2}{2\sigma_{\max}^2}\right) dw \, d\eta_0(z)$$
$$< \infty,$$

where $M_3$ is some positive finite number. The second inequality follows from the inequality $||w| - |\alpha_0 + \alpha_1 z|| \leq |w - \alpha_0 - \alpha_1 z|$.

It is trivial to show that $f^h_{\theta,\eta_1,\eta_2}(y)$ is continuous at $\theta_0$ given $(\eta_1, \eta_2)$ is close to $(\eta_0, \eta_0)$, since $p_\theta(y|z)$ and $h_\theta(y|z)$ are both continuous at $\theta_0$ for $P$-almost every $Y$. Next, we need to show $f^h_{\theta,\eta_1,\eta_2}(y)$ is continuous at $(\eta_0, \eta_0)$ for fixed $\theta$ around $\theta_0$. Accordingly,

$$|f^h_{\theta,\eta_1,\eta_2}(y) - f^h_{\theta,\eta_{10},\eta_{20}}(y)|$$
$$\leq \left|(\eta_1 - \eta_{10})\left(h_\theta \frac{p_\theta}{\eta_2 p_\theta}\right)\right| + \frac{|\eta_{10} h_\theta p_\theta|}{\eta_{20} p_\theta} \left|(\eta_2 - \eta_{20})\left(\frac{p_\theta}{\eta_2 p_\theta}\right)\right|$$
$$\leq \left|(\eta_1 - \eta_{10})\left(\sum_{l=0}^{L} w^l G_l(z;\theta) \frac{p_\theta}{\eta_2 p_\theta}\right)\right|$$
$$+ \sum_{l=0}^{L} \frac{|\eta_{10} g_l(z;\sigma,\gamma,\beta)(w - \alpha_0 - \alpha_1 z)^l p_\theta|}{\eta_{20} p_\theta} \left|(\eta_2 - \eta_{20})\left(\frac{p_\theta}{\eta_2 p_\theta}\right)\right|$$
$$\lesssim \sum_{l=0}^{L} |w|^l \left|(\eta_1 - \eta_{10})\left(G_l(z;\theta) \frac{p_\theta}{\eta_2 p_\theta}\right)\right|$$
$$+ (|w|+1)^l \frac{\eta_{10} p_\theta}{\eta_{20} p_\theta} \left|(\eta_2 - \eta_{20})\left(\frac{p_\theta}{\eta_2 p_\theta}\right)\right|$$
$$\lesssim K_1(w) \times \|\eta_1 - \eta_{10}\|_{BL_1} + K_2(w) \times \|\eta_2 - \eta_{20}\|_{BL_1},$$

where

$$K_1(w) = \exp\left(\frac{M_1|w|}{2\sigma_{\min}^2}\right)\left(\sum_{l=0}^{L}|w|^{l+1} + \sum_{l=0}^{L}|w|^l\right),$$

$$K_2(w) = \exp\left(\frac{M_1|w|}{\sigma_{\min}^2}\right)\sum_{l=0}^{L}(|w|+1)^{l+1}$$

and where $h_\theta$ and $p_\theta$ are abbreviations for $h_\theta(y|z)$ and $p_\theta(y|z)$, respectively. The second inequality follows from

$$h_\theta(y|z) \equiv \sum_{l=0}^{L} g_l(z;\sigma,\gamma,\beta)(w - \alpha_0 - \alpha_1 z)^l = \sum_{l=0}^{L} G_l(z;\theta)w^l,$$

where $G_l(z;\theta) = \sum_{k=l}^{L} g_k(z;\sigma,\gamma,\beta)(-\alpha_0 - \alpha_1 z)^{k-l}(k!/(l!(k-l)!))$. It is trivial to check that $\sum_{i=1}^{N} w_i f_i(z) g_i(z)/2 \in \mathcal{C}_1^1(\mathcal{Z})$ if $f_i(z)$ and $g_i(z)$ belong to $\mathcal{C}_1^1(\mathcal{Z})$. Since the $w_i$'s are nonnegative weights which sum to one, we can find a positive number $R$ such that $R^{-1}G_l(z;\theta) \in \mathcal{C}_1^1(\mathcal{Z})$ for $0 \leq l \leq L$. The last inequality follows from Lemma 6.1 and (43). Note that both $K_1(w)$ and $K_2(w)$ are bounded in $L_1(P)$. This completes the proof. □

VERIFICATION OF ASSUMPTION 4. By repeatedly applying Lemma 6.2, we can check the continuity and boundedness conditions in assumption 4 by resetting $h_\theta(y|z)$ equal to $a_0^T \dot{\ell}_\theta(y|z)$, $a_0^T \ddot{\ell}_\theta(y|z) a_0$ and $a_0^T \dot{\ell}_\theta(y|z) H_0(z)^T a_0$. □

Continuing with the proof of Lemma 6, we need the following verification of assumption 5 (which requires Lemmas 6.3 and 6.4 below).

VERIFICATION OF ASSUMPTION 5. Lemma 6.3 is proved in [27]. The more general version of this lemma can be found on pages 158–159 of [29]. We know the random variable $d$ is binary and thus not smooth. But, if the classes of functions obtained by fixing $d$ to either 0 or 1 are both $P$-Donsker when viewed as functions of the remaining arguments, then the entire classes are $P$-Donsker. A more formal statement of this result can be found in Lemma 9.2 of [23]. Thus, we consider the classes of functions in the following two lemmas for $d = 0$ and $d = 1$, respectively.

LEMMA 6.3. *Let $\mathcal{X} = \bigcup_{j=1}^{\infty} I_j$ be a partition of $\mathbb{R}^1$ into bounded, convex sets whose Lebesgue measure is bounded uniformly away from zero and infinity. Let $\mathcal{G}$ be a class of functions $g : \mathcal{X} \mapsto \mathbb{R}^1$ such that the restrictions $g|_{I_j}$ belong to $C^1_{N_j}$ for every $j$. $\mathcal{G}$ is then $P$-Donsker or $P$-Glivenko–Cantelli for every probability measure $P$ on $\mathcal{X}$ if and only if $\sum_{j=1}^{\infty} N_j P^{1/2}(I_j) < \infty$ or $\sum_{j=1}^{\infty} N_j P(I_j) < \infty$, respectively.*

LEMMA 6.4. *(46) below is P-Donsker when $\theta$ is in some neighborhood of $\theta_0$ and $(\eta_{1j}, \eta_{2j})$ is in some neighborhood of $(\eta_0, \eta_0)$ over compact support $\mathcal{Z}$ for $j = 1, \ldots, k$. The form of $f^{h_j}_{\theta, \eta_{1j}, \eta_{2j}}(y)$ is given in (45) and*

$$(46) \qquad L^H(y; \theta, \bar{\eta}_1, \bar{\eta}_2) \equiv \prod_{j=1}^{k} f^{h_j}_{\theta, \eta_{1j}, \eta_{2j}}(y),$$

*where $H = (h_{\theta 1}(y|z), h_{\theta 2}(y|z), \ldots, h_{\theta k}(y|z))^T$, $\bar{\eta}_1 = (\eta_{11}(z), \ldots, \eta_{1k}(z))^T$, $\bar{\eta}_2 = (\eta_{21}(z), \ldots, \eta_{2k}(z))^T$ and $h_{\theta_j}(y|z) = \sum_{l=0}^{L_j} g_{lj}(z; \sigma, \gamma, \beta)(w - \alpha_0 - \alpha_1 z)^l$, and where $g_{lj}(z; \sigma, \gamma, \beta) \in \mathcal{C}^1_1(\mathcal{Z})$ for $l = 0, 1, \ldots, L_j$ and $j = 1, \ldots, k$.*

PROOF. Without loss of generality, we assume $d = 1$ in the following proof. Based on (43), we have, in each $I_j = \{j - 1 \leq |w| \leq j\}$, $j = 1, \ldots, k$,

$$(47) \qquad |f^{h_j}_{\theta, \eta_{1j}, \eta_{2j}}(y)| \lesssim \sum_{l=0}^{L_j} (|w| + 1)^l \exp\left(\frac{|w|M_1}{2\sigma^2}\right)$$

and

$$(48) \qquad \begin{aligned} &\left|\frac{\partial}{\partial w} f^{h_j}_{\theta, \eta_{1j}, \eta_{2j}}(w, d=1)\right| \\ &\leq \frac{\eta_{1j}(|\frac{\partial}{\partial w} h_{\theta j}|p_\theta)}{\eta_{2j} p_\theta} + \frac{\eta_{1j}(|h_{\theta j} \frac{\partial}{\partial w} \log p_\theta|p_\theta)}{\eta_{2j} p_\theta} \\ &\quad + \frac{\eta_{1j}(|h_{\theta j}|p_\theta)}{\eta_{2j} p_\theta} \times \frac{\eta_{2j}(|\frac{\partial}{\partial w} \log p_\theta|p_\theta)}{\eta_{2j} p_\theta} \\ &\lesssim \sum_{l=0}^{L_j+1} (|w| + 1)^l \left(\exp\left(\frac{|w|M_1}{2\sigma^2}\right) + \exp\left(\frac{|w|M_1}{\sigma^2}\right)\right). \end{aligned}$$

From the above two inequalities, we have that $|(\partial/\partial w) L^H(y; \theta, \bar{\eta}_1, \bar{\eta}_2)|$ is bounded by some constant times $\sum_{l=0}^{R}(j+1)^l(\exp(jM_1 k/2\sigma^2) + \exp(jM_1(k+1)/2\sigma^2))$, where $R = 1 + \sum_{j=1}^{k} L_j$, in each $I_j$, $j \geq 1$. We can then apply Lemma 6.3 to the function $w \mapsto L^H(y; \theta, \bar{\eta}_1, \bar{\eta}_2)$ with $d = 1$ in each $I_j$ defined above. Since the tails in $w$ of $P$ are sub-Gaussian, the series $\sum_j (\sum_{l=0}^{R}(j+1)^l(\exp(jM_1 k/2\sigma^2) + \exp(jM_1(k+1)/2\sigma^2))) P(j-1 \leq |w| \leq j)^{1/2}$ is convergent. Thus, we prove that (46) is $P$-Donsker, which is trivially $P$-Glivenko–Cantelli by Lemma 6.3. □

Continuing with the proof of Lemma 6, we next apply Lemma 6.3 and Lemma 6.4 to show that $x \mapsto \ddot{\ell}(t, \theta, \eta)(x)$ is $P$-Donsker when $(t, \theta, \eta)$ is around $(\beta_0, \theta_0, \eta_0)$. The first term of $\ddot{\ell}(t, \theta, \eta)$, $a_0^T \ddot{\ell}_\theta(y|z) a_0$, is $P$-Donsker, provided the following are both $P$-Donsker, for $0 < r, s, t < \infty$ in (50):

$$f(z): z \mapsto \frac{\exp(\gamma + \beta e^z) e^z}{(1 + \exp(\gamma + \beta e^z))^2}, \tag{49}$$

$$g_{r,s,t}(w): w \mapsto \frac{z^r (w - \alpha_0 - \alpha_1 z)^s}{\sigma^t}. \tag{50}$$

(49) is trivially $P$-Donsker since the function $u \mapsto u \exp(\gamma + \beta u)(1 + \exp(\gamma + \beta u))^{-2}$ is Lipschitz continuous, where $u = e^z$ is $P$-Donsker. For (50), we need to consider Lemma 6.3. We have that $|(\partial/\partial w) g_{r,s,t}(w)| \lesssim (j + |\alpha_0| + |\alpha_1|)^{s-1}$ when $j - 1 \le |w| \le j$. Since the tails in $w$ of $P$ are sub-Gaussian, the series $\sum_j j^{s-1} P(j - 1 \le |w| \le j)^{1/2}$ is convergent. We have thus proven that the first term of $x \mapsto \ddot{\ell}(t, \theta, \eta)(x)$ is $P$-Donsker. By setting $h_\theta(y|z)$ in Lemma 6.4 equal to $a_0^T \dot{\ell}_\theta(y|z)$, $a_0^T \ddot{\ell}_\theta(y|z) a_0$ or $a_0^T H_0(z)$, we can show that the remaining parts of $x \mapsto \ddot{\ell}(t, \theta, \eta)(x)$ are also $P$-Donsker. It can also be proven that $x \mapsto \ell_{t,\theta}(t, \theta, \eta)(x)$ is $P$-Donsker and that $x \mapsto \ell^{(3)}(t, \theta, \eta)(x)$ is $P$-Glivenko–Cantelli by similar reasoning. Thus, assumption 5 is satisfied. The proof is complete since Lemmas 6.5, 6.6 and 6.7 below verify assumption 6. □

LEMMA 6.5. *(2) holds when $\eta$ is in some neighborhood of $\eta_0$.*

PROOF. Based on the form of $\dot{\ell}(\theta, \theta, \eta)$, we can prove (2), provided

$$\mathbb{G}_n(f^h_{\theta,\eta,\eta}(y) - f^h_{\theta,\eta_0,\eta_0}(y)) = O_P(\|\eta - \eta_0\|_{BL_1}). \tag{51}$$

Note that $h_\theta(y|z) = \sum_{l=0}^{L} g_l(z; \sigma, \gamma, \beta)(w - \alpha_0 - \alpha_1 z)^l$ for $\theta \in \Theta$, where $g_l(z; \sigma, \gamma, \beta) \in \mathcal{C}_1^1(\mathcal{Z})$ for $l = 0, 1, \ldots, L$. Thus, (51) will hold, provided

$$\frac{\int g(z) w^l p_{\theta_0}(y|z) \, d\eta(z)}{\int p_{\theta_0}(y|z) \, d\eta(z) \|\eta - \eta_0\|_{BL_1}} - \frac{\int g(z) w^l p_{\theta_0}(y|z) \, d\eta_0(z)}{\int p_{\theta_0}(y|z) \, d\eta_0(z) \|\eta - \eta_0\|_{BL_1}}, \tag{52}$$

for $g(z)$ ranging over $\mathcal{C}_1^1(\mathcal{Z})$, is $P$-Donsker for $l = 0, 1, \ldots, L$. Without loss of generality, it will be enough to verify this for $d = 1$. Note that (52) can also be written as the sum of $Q_{\theta_0, \eta_0, \eta}(w)$ and $-R_{\theta_0, \eta_0, \eta}(w)$, where

$$Q_{\theta_0, \eta_0, \eta}(w) = \frac{\int g(z) w^l p_{\theta_0}(y|z) \, d(\eta - \eta_0)(z)}{\int p_{\theta_0}(y|z) \, d\eta(z) \|\eta - \eta_0\|_{BL_1}}$$

and

$$R_{\theta_0, \eta_0, \eta}(w) = \frac{\int g(z) w^l p_{\theta_0}(y|z) \, d\eta_0(z)}{\int p_{\theta_0}(y|z) \, d\eta_0(z)} \times \frac{\int p_{\theta_0}(y|z) \, d(\eta - \eta_0)(z)}{\int p_{\theta_0}(y|z) \, d\eta(z) \|\eta - \eta_0\|_{BL_1}}.$$

We apply Lemma 6.3 to prove that $Q_{\theta_0,\eta_0,\eta}(w)$ is $P$-Donsker:

$$|(\partial/\partial w)Q_{\theta_0,\eta_0,\eta}(w)|$$

$$\leq l|w|^{l-1}\frac{|(\eta-\eta_0)(gp_{\theta_0})|}{\eta p_{\theta_0}\|\eta-\eta_0\|_{BL_1}} + |w|^l\frac{|(\eta-\eta_0)(gp_{\theta_0}\frac{\partial}{\partial w}\log p_{\theta_0})|}{\eta p_{\theta_0}\|\eta-\eta_0\|_{BL_1}}$$

$$+ |w|^l\frac{|(\eta-\eta_0)(gp_{\theta_0})|}{\eta p_{\theta_0}\|\eta-\eta_0\|_{BL_1}} \times \frac{\eta(p_{\theta_0}|\frac{\partial}{\partial w}\log p_{\theta_0}|)}{\eta p_{\theta_0}}$$

$$\lesssim \sum_{m=l-1}^{l+1}|w|^m\frac{|(\eta-\eta_0)(s_m(z)p_{\theta_0}/(\eta p_{\theta_0}))|}{\|\eta-\eta_0\|_{BL_1}},$$

where $g$ and $p_{\theta_0}$ are abbreviations for $g(z)$ and $p_{\theta_0}(y|z)$, respectively, and $s_m(z) \in \mathcal{C}_1^1(\mathcal{Z})$ for $m = l-1, l, l+1$. Combining this with (43), $|(\partial/\partial w)Q_{\theta_0,\eta_0,\eta}(w)|$ is bounded by a constant times $\exp(jM_1/2\sigma^2)\sum_{m=l-1}^{l+2}j^m$ in each region $I_j = \{j-1 \leq |w| \leq j\}$. It is thus proved that $Q_{\theta_0,\eta_0,\eta}(w)$ is $P$-Donsker, by Lemma 6.3. Similarly, we can also show that $R_{\theta_0,\eta_0,\eta}(w)$ is $P$-Donsker. This completes the proof. □

LEMMA 6.6. *(3) and (4) hold when $\eta$ is in some neighborhood of $\eta_0$.*

PROOF. Based on the form of $\ddot{\ell}(\theta,\theta,\eta)$, (3) will follow provided

$$(53) \quad P\left|\frac{\eta(g(z)w^l p_{\theta_0})}{\eta p_{\theta_0}} - \frac{\eta_0(g(z)w^l p_{\theta_0})}{\eta_0 p_{\theta_0}}\right| = O_P(\|\eta-\eta_0\|_{BL_1}),$$

for any $g(z) \in \mathcal{C}_1^1(\mathcal{Z})$ and for $l = 0, 1, \ldots, L$. Now, (53) is bounded by the summation of $P|\bar{Q}_{\theta_0,\eta_0,\eta}(w)|$ and $P|\bar{R}_{\theta_0,\eta_0,\eta}(w)|$, where $\bar{Q}_{\theta_0,\eta_0,\eta}(w) \equiv Q_{\theta_0,\eta_0,\eta}(w)\|\eta-\eta_0\|_{BL_1}$ and $\bar{R}_{\theta_0,\eta_0,\eta}(w) \equiv R_{\theta_0,\eta_0,\eta}(w)\|\eta-\eta_0\|_{BL_1}$, and where $Q_{\theta_0,\eta_0,\eta}(w)$ and $R_{\theta_0,\eta_0,\eta}(w)$ are as defined in the proof of Lemma 6.5 above. Note that $P|\bar{Q}_{\theta_0,\eta_0,\eta}(w)|$ can be written as

$$P|\bar{Q}_{\theta_0,\eta_0,\eta}(w)|$$

$$= \int_\mathbb{R} \Delta_0(w)|w|^l\left|\int_\mathcal{Z} g(z)p_{\theta_0}(w,d=0|z)\,d(\eta-\eta_0)(z)\right| dw\, P(d=0)$$

$$+ \int_\mathbb{R} \Delta_1(w)|w|^l\left|\int_\mathcal{Z} g(z)p_{\theta_0}(w,d=1|z)\,d(\eta-\eta_0)(z)\right| dw\, P(d=1),$$

where $\Delta_i(w) = \int p_{\theta_0}(w,d=i|z)\,d\eta_0(z)/\int p_{\theta_0}(w,d=i|z)\,d\eta(z)$ for $i = 0, 1$. Without loss of generality, we can show $P|\bar{Q}_{\theta_0,\eta_0,\eta}(w)|$ is of the order $\|\eta-\eta_0\|_{BL_1}$, provided the first integral on the right-hand side of the above equation is of the same order. Based on the inequality $||w|-|\alpha_0+\alpha_1 z|| \leq |w-\alpha_0-\alpha_1 z| \leq |w|+|\alpha_0+\alpha_1 z|$, we have $\Delta_0(w) \lesssim \exp((M_1/2\sigma^2)|w|)$. We can

verify that $\exp(w^2/4\sigma^2)g(z)p_{\theta_0}(w,d=0|z) \in C_1^1(Z)$ for $P$-almost all $Y$ after proper rescaling. Also, $\int_{\mathbb{R}} \Delta_0(w)|w|^l \exp(-w^2/4\sigma^2)\,dw$ is trivially bounded. Thus, $P|\bar{Q}_{\theta_0,\eta_0,\eta}(w)|$ is of the same order as $\|\eta - \eta_0\|_{BL_1}$. By similar analysis, we can also show that $P|\bar{R}_{\theta_0,\eta_0,\eta}(w)| = O_P(\|\eta - \eta_0\|)$. Since $\ell_{t,\theta}(t,\theta,\eta)$ is similar to $\ddot{\ell}(t,\theta,\eta)$, we also have $P\ell_{t,\theta}(\theta_0,\theta_0,\eta) - P\ell_{t,\theta}(\theta_0,\theta_0,\eta_0) = O_P(\|\eta - \eta_0\|_{BL_1})$. □

LEMMA 6.7. *(5) holds when $\eta$ is in some neighborhood of $\eta_0$.*

PROOF. Based on previous discussions about the verification of (5), we only need to show that $|\dot{\ell}(\theta_0,\theta_0,\eta) - \dot{\ell}(\theta_0,\theta_0,\eta_0)| \lesssim C(y)\|\eta - \eta_0\|_{BL_1}$ and $|\text{lik}(\theta_0,\eta) - \text{lik}(\theta_0,\eta_0) - A_0(\eta - \eta_0)\text{lik}(\theta_0,\eta_0)| \lesssim D(y)\|\eta - \eta_0\|_{BL_1}^2$, where $C(y)$ and $D(y)$ are both bounded in $L_2(P)$. The former inequality is easily proved via techniques similar to those used in the proof of Lemma 6.6. For the latter, we can write

$$\text{lik}(\theta_0,\eta) - \text{lik}(\theta_0,\eta_0)$$
$$= A_0(\eta - \eta_0)\text{lik}(\theta_0,\eta_0) + p_{\theta_0}(y|z)\int_{\{z_c\}} d(\eta - \eta_0)\int p_{\theta_0}(y|z)\,d(\eta - \eta_0),$$

where $A_0(\eta - \eta_0) = \eta_0\{z_c\}^{-1}\int_{\{z_c\}} d(\eta - \eta_0) + (\int p_{\theta_0}(y|z)\,d\eta_0(z))^{-1} \times \int p_{\theta_0}(y|z)\,d(\eta - \eta_0)$. It is now easy to show that $|p_{\theta_0}(y|z)\int 1_{\{z_c\}}\,d(\eta - \eta_0) \times \int p_{\theta_0}(y|z)\,d(\eta - \eta_0)| \lesssim D(y)\|\eta - \eta_0\|_{BL_1}^2$ since $p_{\theta_0}(y|z) \in \mathcal{C}_1^1(\mathcal{Z})$ for $P$-almost every $Y$ via rescaling. □

PROOF OF LEMMA 7. The proof is analogous to that of Lemma 3 in [17]. □

**Acknowledgments.** The authors thank the referees and editors for several helpful suggestions.

## REFERENCES

[1] BEGUN, J. M., HALL, W. J., HUANG, W.-.M. and WELLNER, J. A. (1983). Information and asymptotic efficiency in parametric–nonparametric models. *Ann. Statist.* **11** 432–452. MR0696057
[2] BICKEL, P. J., KLAASSEN, C. A. J., RITOV, Y. and WELLNER, J. A. (1998). *Efficient and Adaptive Estimation for Semiparametric Models.* Springer, New York. MR1623559
[3] CHENG, G. (2006). Higher order semiparametric frequentist inference and the profile sampler. Ph.D. thesis, Dept. Statistics, Univ. Wisconsin-Madison.
[4] COX, D. R. (1972). Regression models and life-tables. *J. Roy. Statist. Soc. Ser. B* **34** 187–220. MR0341758

[5] Dalalyan, A. S., Golubev, G. K. and Tsybakov, A. B. (2006). Penalized maximum likelihood and semiparametric second-order efficiency. *Ann. Statist.* **34** 169–201. MR2275239

[6] Ghosh, J. K. and Mukerjee, R. (1991). Characterization of priors under which Bayesian and frequentist Bartlett corrections are equivalent in the multiparameter case. *J. Multivariate Anal.* **38** 385–393. MR1131727

[7] Groeneboom, P. (1991). Nonparametric maximum likelihood estimators for interval censoring and deconvolution. Technical Report No. 378, Dept. Statistics, Stanford Univ.

[8] Hall, P. (1992). *The Bootstrap and Edgeworth Expansion*. Springer, New York. MR1145237

[9] Härdle, W. and Tsybakov, A. B. (1993). How sensitive are average derivatives? *J. Econometrics* **58** 31–48. MR1230979

[10] Huang, J. (1996). Efficient estimation for the Cox model with interval censoring. *Ann. Statist.* **24** 540–568. MR1394975

[11] Johnson, R. A. (1970). Asymptotic expansions associated with posterior distributions. *Ann. Math. Statist.* **41** 851–864. MR0263198

[12] Kalbfleisch, J. D. (1978). Nonparametric Bayesian analysis of survival time data. *J. Roy. Statist. Soc. Ser. B* **40** 214–221. MR0517442

[13] Kosorok, M. R. (2008). *Introduction to Empirical Processes and Semiparametric Inference*. Springer, New York.

[14] Kosorok, M. R., Lee, B. L. and Fine, J. P. (2004). Robust inference for univariate proportional hazards frailty regression models. *Ann. Statist.* **32** 1448–1491. MR2089130

[15] Kuo, H. H. (1975). *Gaussian Measure on Banach Spaces. Lecture Notes in Math.* **463**. Springer, Berlin. MR0461643

[16] Lawley, D. N. (1956). A general method for approximating to the distribution of likelihood ratio criteria. *Biometrika* **43** 295–303. MR0082237

[17] Lee, B. L., Kosorok, M. R. and Fine, J. P. (2005). The profile sampler. *J. Amer. Statist. Assoc.* **100** 960–969. MR2201022

[18] Ma, S. and Kosorok, M. R. (2005). Robust semiparametric M-estimation and the weighted bootstrap. *J. Multivariate Analysis* **96** 190–217. MR2202406

[19] Murphy, S. A., Rossini, A. J. and van der Vaart, A. W. (1997). MLE in the proportional odds model. *J. Amer. Statist. Assoc.* **92** 968–976. MR1482127

[20] Murphy, S. A. and van der Vaart, A. W. (1997). Semiparametric likelihood ratio inference. *Ann. Statist.* **25** 1471–1509. MR1463562

[21] Murphy, S. A. and Van der Vaart, A. W. (1999). Observed information in semiparametric models. *Bernoulli* **5** 381–412. MR1693616

[22] Murphy, S. A. and Van der Vaart, A. W. (2000). On profile likelihood (with discussion). *J. Amer. Statist. Assoc.* **95** 449–485. MR1803168

[23] Murphy, S. A. and Van der Vaart, A. W. (2001). Semiparametric mixtures in case-control studies. *J. Multivariate Anal.* **79** 1–32. MR1867252

[24] Roeder, K., Carroll, R. J. and Lindsay, B. G. (1996). A semiparametric mixture approach to case-control studies with errors in covariables. *J. Amer. Statist. Assoc.* **91** 722–732. MR1395739

[25] Severini, T. A. and Wong, W. H. (1992). Profile likelihood and conditionally parametric models. *Ann. Statist.* **20** 1768–1802. MR1193312

[26] Shen, X. (2002). Asymptotic normality in semiparametric and nonparametric posterior distributions. *J. Amer. Statist. Assoc.* **97** 222–235. MR1947282

[27] VAN DER VAART, A.W. (1996). New Donsker classes. *Ann. Probab.* **24** 2128–2140. MR1415244
[28] VAN DER VAART, A. W. (1998). *Asymptotic Statistics*. Cambridge Univ. Press. MR1652247
[29] VAN DER VAART, A. W. and WELLNER, J. A. (1996). *Weak Convergence and Empirical Processes. With Applications to Statistics*. Springer, New York. MR1385671
[30] WELCH, B. L. and PEERS, H. W. (1963). On formulae for confidence points based on integrals of weighted likelihood. *J. Roy. Statist. Soc. Ser. B* **25** 318–329. MR0173309
[31] WELLNER, J. A. and ZHANG, Y. (2007). Two likelihood-based semiparametric estimation methods for panel count data with covariates. *Ann. Statist.* **35** 2106–2142. MR2363965
[32] WOUK, A. (1979). *A Course of Applied Functional Analysis*. Wiley, New York. MR0524723

INSTITUTE OF STATISTICS
AND DECISION SCIENCES
DUKE UNIVERSITY
214 OLD CHEMISTRY BUILDING
DURHAM, NORTH CAROLINA 27708
USA
E-MAIL: chengg@stat.duke.edu

DEPARTMENT OF BIOSTATISTICS
UNIVERSITY OF NORTH CAROLINA
AT CHAPEL HILL
3101 MCGAVRAN-GREENBERG HALL
CHAPEL HILL, NORTH CAROLINA 27599
USA
E-MAIL: kosorok@unc.edu